\documentstyle[12pt,amssymb]{amsart}

\numberwithin{equation}{section}

\newtheorem{th}{Th\'eor\`eme}[section]
\newtheorem{cor}[th]{Corollaire}
\newtheorem{lemme}[th]{Lemme}
\newtheorem{prop}[th]{Proposition}

\theoremstyle{definition}
\newtheorem{defi}{D\'efinition}[section]
\newtheorem{question}[th]{Question}

\mathcode`A="7041 \mathcode`B="7042 \mathcode`C="7043 \mathcode`D="7044
\mathcode`E="7045 \mathcode`F="7046 \mathcode`G="7047 \mathcode`H="7048
\mathcode`I="7049 \mathcode`J="704A \mathcode`K="704B \mathcode`L="704C
\mathcode`M="704D \mathcode`N="704E \mathcode`O="704F \mathcode`P="7050
\mathcode`Q="7051 \mathcode`R="7052 \mathcode`S="7053 \mathcode`T="7054
\mathcode`U="7055 \mathcode`V="7056 \mathcode`W="7057 \mathcode`X="7058
\mathcode`Y="7059 \mathcode`Z="705A

\mathchardef\nearrow"3225
\mathchardef\searrow"3226
\mathchardef\nwarrow"322D
\mathchardef\swarrow"322E

\newcommand{\hfl}[2]{\smash{\mathop{\hbox to 6mm{\rightarrowfill}}
                     \limits^{\scriptstyle#1}_{\scriptstyle#2}}}

\newcommand{\vfl}[2]{\llap{$\scriptstyle #1$}\left\downarrow
                     \vbox to 3mm{}\right.\rlap{$\scriptstyle #2$}}

\newcommand{\efl}[2]{\smash{\mathop{\hbox to 6mm{\rightarrowfill}}
                     \limits^{\scriptstyle#1}_{\scriptstyle#2}}}

\newcommand{\sfl}[2]{\llap{$\scriptstyle #1$}\left\downarrow
                     \vbox to 3mm{}\right.\rlap{$\scriptstyle #2$}}

\newcommand{\nfl}[2]{\llap{$\scriptstyle #1$}\left\uparrow
                     \vbox to 3mm{}\right.\rlap{$\scriptstyle #2$}}

\newcommand{\wfl}[2]{\smash{\mathop{\hbox to 6mm{\leftarrowfill}}
                     \limits^{\scriptstyle#1}_{\scriptstyle#2}}}

\newcommand{\nefl}[2]{\llap{$\scriptstyle#1$}\hspace{-2pt}%
                  \nearrow\hspace{-2pt}\rlap{$\scriptstyle #2$}}

\newcommand{\sefl}[2]{\raise2pt\hbox{\llap{$\scriptstyle#1$}}\hspace{-2pt}%
                  \searrow\hspace{-2pt}\raise2pt\hbox{\rlap{$\scriptstyle#2$}}}

\newcommand{\swfl}[2]{\raise2pt\hbox{\llap{$\scriptstyle#1$}}\hspace{-2pt}
                  \swarrow\hspace{-2pt}\raise2pt\hbox{\rlap{$\scriptstyle#2$}}}

\newenvironment{diagram}{\begin{matrix}}{\end{matrix}}

\newcommand{\nc}{\newcommand}
\nc{\on}{\operatorname}
\nc{\onw}{\operatornamewithlimits}

\nc{\reln}{{\Bbb Z}}
\nc{\quot}{{\Bbb Q}}
\nc{\comp}{{\Bbb C}}
\nc{\proj}{{\Bbb P}}
\nc{\real}{{\Bbb R}}
\nc{\natn}{{\Bbb N}}
\nc{\planp}{\proj^2}

\nc{\tr}{\Bbb}

\nc{\isom}{\cong}

\nc{\lra}{\longrightarrow}
\nc{\ra}{\rightarrow}
\nc{\GL}{\on{GL}}
\nc{\SL}{\on{SL}}
\nc{\U}{\on{U}}
\nc{\SO}{\on{SO}}
\nc{\Spin}{\on{Spin}}

\nc{\Sym}{\on{Sym}}
\nc{\Res}{\onw{Res}}
\nc{\Ext}{\on{Ext}}
\nc{\Tor}{\on{Tor}}
\nc{\Hom}{\on{Hom}}
\nc{\End}{\on{End}}
\nc{\Pic}{\on{Pic}}
\renewcommand{\H}{\on{H}}
\nc{\Todd}{\on{Todd}}
\nc{\res}{\on{Res}}
\nc{\Spec}{\on{Spec}}
\nc{\Gr}{\on{Gr}}

\nc{\Tr}{\on{Tr}}
\nc{\tensor}{\otimes}
\nc{\rang}{\on{rang}}
\nc{\ds}{\displaystyle}
\nc{\supp}{\on{supp}}
\nc{\supps}{\supp_{s}}
\nc{\prof}{\on{prof}}
\nc{\dh}{\on{dh}}
\nc{\gr}{\on{gr}}
\nc{\id}{\on{id}}

\nc{\HyperExt}{\on{\tr{E}}xt}
\nc{\HyperH}{\on{\tr{H}}}
\renewcommand{\Im}{\on{Im}}
\nc{\Ker}{\on{Ker}}
\nc{\Coker}{\on{Coker}}
\nc{\Aut}{\on{Aut}}
\nc{\codim}{codim}

\nc{\Cl}{\on{Cl}}
\nc{\Pf}{\on{Pf}}
\nc{\ext}{\on{ext}}
\nc{\h}{\on{h}}
\nc{\mult}{\on{mult}}
\nc{\Groth}{\on{Groth}}
\nc{\Stab}{\on{Stab}}
\nc{\droitep}{\proj_{1}}

\newbox\boxa
\nc{\boxit}[1]{\setbox\boxa=\hbox{#1}%
  \setbox\boxa=\hbox{\vrule\box\boxa\vrule}%
  \setbox\boxa=\vbox{\hrule\box\boxa\hrule}%
  \box\boxa\relax}
\nc{\etimes}{\mathop{\raise -1pt\hbox{\boxit{$\times$}}}}

\nc{\ie}{{\em ie. }}
\nc{\ul}{\underline}
\nc{\ol}{\overline}
\newenvironment{proof}{\begin{pf}}{\end{pf}}
\newcounter{lc}
\renewcommand{\restriction}[1]{\,\raisebox{-1ex}{$\mid_{#1}$}}
\nc{\inject}{\hookrightarrow}

\nc{\cf}{{\em cf.} }
\nc{\ort}{\perp}
\nc{\R}{\on{R}}
\nc{\cqfd}{\qed}
\nc{\tvi}{\vrule height 12pt depth 5pt width 0pt}
\nc{\tv}{\tvi\vrule}
\nc{\osum}{\onw{\oplus}}

\nc{\np}{\clearpage}

\renewcommand{\codim}{\on{codim}}

\title[Groupe de Picard]{Groupe de Picard de la variete des
theta-caracteristiques
des courbes planes}
\author{Christoph Sorger}
\address{Institut de math\'ematiques Jussieu\\Universit\'e Denis Diderot - Case
Postale 7012\\2, place Jussieu\\75251 Paris CEDEX 05}
\email{sorger@@mathp7.jussieu.fr}
\thanks{Pr\'e-publication les pour annales de la conf\'erence annuelle 93
 d'Europroj \'a Catania}

\begin{document}
\addtolength{\baselineskip}{2pt}
\maketitle
\begin{abstract}
Nous calculons le groupe de Picard de l'espace des modules
$\Theta_{\planp}(d)$ des th\^eta-caract\'eristiques des courbes planes (non
forc\'ement
lisses) de degr\'e $d$. Nous montrons que pour $d\ge 6$, le groupe de Picard de
la
composante paire est engendr\'e par le fibr\'e d\'eterminant ${\cal{D}}$ et
l'image
r\'eciproque ${\cal{L}}$, sous le morphisme support sch\'ematique, de
${\cal{O}}_{\proj}(1)$ o\`u $\proj$ est l'espace des courbes de degr\'e $d$ sur
$\planp$.  Pour $1\le d \le 4$, ce groupe est engendr\'e uniquement par
${\cal{L}}$,
le fibr\'e d\'eterminant \'etant trivial dans ces cas.  Dans le groupe $\Cl$
des diviseurs
de Weil, le fibr\'e d\'eterminant admet une racine: le fibr\'e pfaffien
${\cal{P}}$. De
plus, pour $d\geq 4$ pair,
${\cal{L}}$ aussi admet une racine dans $\Cl$. Si $d=5$, le fibr\'e pfaffien
s'\'etend
en un fibr\'e inversible et le groupe de Picard est engendr\'e par ${\cal{L}}$
et
${\cal{P}}$. On verra ainsi que la composante paire de
$\Theta_{\planp}(d)$ est localement factoriel si et seulement si $d=1,2,3$ ou
$d=5$.
On obtient un r\'esultat analogue pour la composante impaire.
Ensuite, nous \'etudions la question d'existence d'une famille universelle sur
un
ouvert de l'ouvert $U$ des th\^eta-caract\'eristiques dont le faisceau
sous-jacent est
stable. Nous montrons que pour $d$ pair une telle famille ne peut exister,
tandis
que pour $d$ impair une telle famille existe, pas sur $U$ entier, mais
localement
dans la topologie de Zariski.
\end{abstract}

\section{Introduction}
Le pr\'esent travail poursuit l'\'etude de l'espace des modules de \cite{11}
des
th\^eta-caract\'eristiques des courbes (non forc\'ement lisses) trac\'ees sur
une surface
dans le cas du plan projectif. Nous calculons le groupe de Picard
des composantes irr\'eductibles et \'etudions les questions de
factorialit\'e locale et de l'existence d'une famille universelle.

Rappelons qu'on appelle th\^eta-caract\'eristique (g\'en\'eralis\'ee) d'une
courbe plane un ${\cal{O}}_{\planp}$-module ${\cal{F}}$ de
dimension $1$ muni d'une identification sym\'etrique $\sigma$
avec son
$\omega_{_{\planp}}$-dual. La courbe associ\'ee \`a $({\cal{F}},\sigma)$
est d\'efinie par le support sch\'ematique de ${\cal{F}}$, \ie le
$0$-\`eme id\'eal de Fitting de ${\cal{F}}$. Le degr\'e de cette courbe
s'appelle multiplicit\'e de ${\cal{F}}$. D'apr\`es $\cite{11}$,
il existe un espace  de modules grossier $\Theta_{\planp}(d)$ pour
les th\^eta-caract\'eristiques semi-stables (\cf \ref{theta-def} pour
la d\'efinition) de multiplicit\'e $d$.
C'est une vari\'et\'e normale et projective de dimension
$\frac{d(d+3)}{2}$, muni d'un morphisme
$$\sigma_{_{\Theta}}:\Theta_{\planp}(d)\lra\proj^{M}$$ avec
$M=\frac{d(d+3)}{2}$, d\'efini
en associant \`a la classe d'une th\^eta-caract\'eristique semi-stable
son support sch\'ematique. Au-dessus de l'ouvert des
courbes lisses, $\sigma_{_{\Theta}}$ est \'etale de degr\'e $2^{2g}$,
la fibre au dessus d'une courbe lisse $C$ s'identifiant \`a l'ensemble
de th\^eta-caract\'eristiques classiques sur $C$.
Par contre, ce morphisme n'est pas quasi-fini: si $C$ est int\`egre
la fibre est finie si et seulement si $C$ n'a que des singularit\'es
simples (\cf Cook dans le m\^eme volume); si $C$ est \`a structure
multiple d'une courbe lisse (\ie $C=rC'$ en tant que diviseur de
Weil avec $C'$ lisse), la fibre contient l'espace des modules de
fibr\'es $\omega_{C'}$-quadratiques de rang $r$ sur la courbe $C'$.

Si $({\cal{F}},\sigma)$ est une famille de th\^eta-caract\'eristiques
param\'etr\'ee par une vari\'et\'e connexe, la dimension de l'espace des
sections de ${\cal{F}}_{s}$ est invariant modulo 2 pour $s\in S$.
(\cite{11}, 0.3 ou  \cite{12}). Il s'ensuit que
$\Theta_{\planp}(d)$ a au moins deux composantes (qui seront non
vide pour $d\ge 4$) correspondant aux th\^eta-caract\'eristiques
$({\cal{F}},\sigma)$  telles que
$\h^{0}({\cal{F}})=0
\bmod 2$ et aux th\^eta-caract\'eristiques  telles que
$\h^{0}({\cal{F}})=1 \bmod 2$. Si $d$ est impair on a une troisi\`eme
composante, qu'on note
$\Theta_{can}(d)$, correspondant aux th\^eta-caract\'eristiques {\em
canoniques},
\ie de la forme ${\cal{O}}_{C}(\frac{d-3}{2})$, o\`u $C$ est une
courbe de degr\'e $d$. Cette composante est toujours isomorphe \`a
l'espace des courbes de degr\'e
$d$ sur $\planp$. Dans ce qui suit, on appellera {\em paire} une
th\^eta-caract\'eristique  non canonique telle que $\h^{0}({\cal{F}})=0
\bmod 2$. Une th\^eta-caract\'eristique non canonique telle que
$\h^{0}({\cal{F}})=1
\bmod 2$ sera appel\'e {\em impaire}. On notera $\Theta_{p}(d)$ (resp.
$\Theta_{i}(d)$) la composante des th\^eta-caract\'eristiques paires
(resp. impaires). D'apr\`es  (\cite{11}, 0.5; \cite{1}
pour l'ouvert des courbes lisses) ces composantes sont
irr\'eductibles.  La composante paire contient l'ouvert not\'e $I(d)$
des th\^eta-caract\'eristiques dites {\em ineffectives} (\ie
$\h^{0}({\cal{F}})=0$), la composante impaire contient l'ouvert
not\'e $SI(d)$ des th\^eta-caract\'eristiques dites {\em
semi-ineffectives} (\ie
$\h^{0}({\cal{F}})=1$).

Si $d=1$, toute th\^eta-caract\'eristique est canonique; si $d=2$,
$\Theta_{p}(d)$ est isomorphe \`a
$\proj_{5}$ et $\Theta_{i}(d)$ est vide. Enfin, si
$d=3$, $\Theta_{i}(d)$ est encore vide, toute th\^eta-caract\'eristique
telle que $\h^{0}({\cal{F}})$ soit impaire \'etant
canonique.

Par normalit\'e de $\Theta_{p}(d)$, le morphisme naturel
$$\rho:\Pic(\Theta_{p}(d))\lra\Cl(\Theta_{p}(d))$$
est une injection.

\begin{th}\label{PicPair} Si $d=3,4$, le
groupe de Picard de $\Theta_{p}(d)$ est un groupe ab\'elien libre \`a
un g\'en\'erateur. De plus, $\rho$ est
un isomorphisme pour $d=3$; si $d=4$, il est d'indice $2$.
Si $d\geq 5$, le groupe de Picard de $\Theta_{p}(d)$ est isomorphe \`a un
groupe ab\'elien libre \`a deux g\'en\'erateurs. De plus, $\rho$
est un isomorphisme pour $d=5$, d'indice $4$ si $d\geq 6$ est pair et
d'indice $2$ si $d\geq 7$ est impair.\par
En particulier, la vari\'et\'e
$\Theta_{p}(d)$ est localement factorielle si et seulement si $d=1,2,3$ ou
$d=5$.
\end{th}

Pour l'identification des g\'en\'erateurs on consid\'era l'espace de
modules $N_{\planp}(d,\chi)$ des ${\cal{O}}_{\planp}$-modules
semi-stables de dimension $1$, de degr\'e $d$ et de caract\'eristique
d'Euler-Poincar\'e $\chi$. Consid\'erons le morphisme d'oubli
$$\beta:\Theta_{\planp}(d)\lra N_{\planp}(d,0)$$ qui associe \`a la
classe d'une th\^eta-caract\'eristique semi-stables la classe du
faisceau semi-stable sous-jacent.
D'apr\`es \cite{9}, $N_{\planp}(d,\chi)$ est une vari\'et\'e
irr\'eductible, projective et normale. En outre,
$N_{\planp}(d,\chi)$ est localement factorielle et son groupe de
Picard est un groupe ab\'elien libre \`a deux g\'en\'erateurs, not\'es
${\cal{L}}_{N}$ et
${\cal{D}}_{N}$, d\'efinies fonctoriellement. Le fibr\'e inversible
${\cal{L}}_{N}$ s'identifie \`a l'image r\'eciproque de
${\cal{O}}_{\proj^{M}}(1)$ sous le morphisme support sch\'ematique
$\sigma_{_{M}}:N_{\planp}(d,0)\lra\proj^{M}$; le fibr\'e
${\cal{D}}_{N}$ est le fibr\'e d\'eterminant (\cf \ref{LePotierspace} pour la
d\'efinition) Notons ${\cal{L}}$ et ${\cal{D}}$ les images r\'eciproques de
${\cal{L}}_{N}$  et ${\cal{D}}_{N}$ sous $\beta$. Pour $d=5$,
le fibr\'e inversible ${\cal{D}}$ admet une racine:
le  fibr\'e pfaffien, not\'e ${\cal{P}}$ (\cf \ref{defpfaff} pour la
d\'efinition). Alors on a pour le groupe de Picard
$$\Pic(\Theta_{p}(d))=
\begin{cases}
	<{\cal{L}}>& \text{si $d=3,4$}\\
 <{\cal{L}},{\cal{P}}>& \text{si $d=5$}\\
 <{\cal{L}},{\cal{D}}>& \text{si $d\ge6$}\\
\end{cases}
$$

Soit $\Theta^{os}_{p}\subset\Theta_{p}$ l'ouvert des
 th\^eta-caract\'eristiques $({\cal{F}},\sigma)$ tels que
${\cal{F}}$ soit stable en tant que faisceau. D'apr\`es \cite{11},
cet ouvert est lisse et l'on a
$$\Cl(\Theta_{p})=\Pic(\Theta_{p}^{os}(d)).$$
On verra aussi que si $d\geq 4$ est pair, l'image de
${\cal{L}}$ sous $\rho$ admet une racine dans $\Pic(\Theta_{p}^{os}(d))$,
not\'e ${\cal{R}}$. Si $d\geq 5$, l'image de ${\cal{D}}$ sous $\rho$
admet une racine dans $\Pic(\Theta_{p}^{os}(d))$,
le fibr\'e pfaffien ${\cal{P}}$. On verra qu'en fait
$$\Pic(\Theta_{p}^{os}(d))=
\begin{cases}
	<{\cal{L}}>& \text{ si $d=3$}\\
	<{\cal{R}}>& \text{ si $d=4$}\\
	<{\cal{L}},{\cal{P}}>& \text{ si $d\geq 5$ est impair}\\
 <{\cal{R}},{\cal{P}}>& \text{ si $d\geq 6$ est pair}\\
\end{cases}
$$
Ceci pr\'ecise l'assertion concernant la factorialit\'e locale.

Pour la composante impaire, on obtient

\begin{th}\label{PicImpair}\par Si $d\ge4$, avec $d\not=6$, le
groupe de Picard de
$\Theta_{i}(d)$ s'identifie \`a un groupe ab\'elien libre \`a deux
g\'en\'erateurs: les images r\'eciproques de ${\cal{L}}_{N}$ et
${\cal{D}}_{N}$. De plus, le morphisme naturel
$\Pic(\Theta_{i}(d))\lra\Cl(\Theta_{i}(d))$ est d'indice $2$.
Si $d=6$, le groupe de
Picard de $\Theta_{i}(d)$ a un g\'en\'erateur de plus provenant du
ferm\'e des th\^eta-caract\'eristiques telles que
$h^{0}({\cal{F}})\ge3$, qui est un diviseur irr\'eductible pour cette
valeur de
$d$.\par En particulier, pour $d\geq 4$, la vari\'et\'e $\Theta_{i}(d)$ n'est
pas localement factorielle.
\end{th}

Nous d\'emontrerons les th\'eor\`emes ci-dessus en \'etudiant
\begin{list}{-}{}
\item le groupe de Picard des ouverts $I(d)\subset\Theta_{p}(d)$
et $SI(d)\subset\Theta_{i}(d)$
\item les sous-sch\'emas correspondant aux th\^eta-caract\'eristiques
sp\'eciales (\ie ni ineffectives, ni semi-ineffectives).
\end{list}

Le groupe de Picard de $I(d)$ et $SI(d)$ se calculera en utilisant
le type particulier de la r\'esolution minimale des
th\^eta-caract\'eristique ineffectives et semi-ineffectives. Si le
support sch\'ematique est lisse, ce type particulier
est connu de Dixon \cite{3}, Catanese \cite{2} et
Laszlo \cite{7}. La d\'emonstration de Laszlo se g\'en\'eralise
facilement au cas des th\^eta-caract\'eristiques \`a support sch\'ematique
quelconque. Cependant, nous en donnons encore une autre
d\'emonstration permettant d'obtenir un r\'esultat davantage fonctoriel
et une description de la fl\`eche en termes de la cohomologie d'une
th\^eta-caract\'eristique donn\'ee, n\'ecessaire pour notre propos.

Pour l'\'etude des th\^eta-caract\'eristiques sp\'eciales nous d\'efinissons
les ensembles
$$\Theta^{r}(d)=\{[{\cal{F}},\sigma]\in\Theta(d)/\h^{0}({\cal{F}})\ge
r+1,
\h^{0}({\cal{F}})=(r+1)\bmod 2\}$$ qu'on munit de leur structure
sch\'ematique de vari\'et\'e d\'eterminantielle naturelle
d\'efinie dans la section \ref{BNlocus}. On y re-d\'emontre aussi que
$$\text{codim}_{\Theta}(\Theta^{r}(d))\le\frac{r(r+1)}{2},$$
si $\Theta^{r}(d)$ est non vide. En particulier $\Theta^{1}(d)$ est un
diviseur, si
non vide, puisque ce n'est pas une composante. L'\'enonc\'e de type
Brill-Noether
duquel on a besoin est le suivant:

\begin{th}\label{BN} Le diviseur $\Theta^{1}(d)$ est irr\'eductible et le
ferm\'e $\Theta^{3}(d)$ est de codimension au moins deux pour tout
$d$. Le ferm\'e $\Theta^{2}(d)$ est de codimension au moins deux sauf
si $d=6$ o\`u il est irr\'eductible de codimension un.
\end{th}

Nous d\'emontrons ce th\'eor\`eme en \'etudiant les ouverts $U_{a}(d)$ de
$\Theta(d)$ des th\^eta-caract\'eristiques de passant pas par un point donn\'e
$a\in\planp$. Le point essentiel ici est le fait que ces ouverts
ont une description particuli\`erement simple en
termes de fibr\'es de Higgs quadratiques sur la droite projective

Enfin, nous \'etudierons la question de l'existence d'une famille
universelle: Soit $U$ un ouvert non-vide de l'ouvert des
th\^eta-caract\'eristiques
${\cal{O}}$-stables paires (resp. impairs). On appelle famille
universelle sur $U$ une famille $({\cal{F}},\sigma)$ param\'etr\'ee par
$U$ telle que le morphisme modulaire induit $U\lra\Theta^{os}(d)$ soit
l'inclusion.

\begin{th}\label{UniFam} Soit $d\geq 3$. Si $d$ est pair il n'existe
pas de famille universelle sur $U$. Si $d$ est impair il existe une
famille universelle localement dans la topologie de Zariski, mais il
n'existe pas de famille universelle globalement sur
$\Theta^{os}(d)$.
\end{th}

{\em Notations:} Par vari\'et\'e alg\'ebrique on entend sch\'ema de type fini,
s\'epar\'e, sur un corps $k$
suppos\'e alg\'ebriquement clos de caract\'eristique 0.

Si
$f:X\lra S$ est un morphisme de vari\'et\'es alg\'ebriques, on note $X_{s}$ la
fibre de $f$ au-dessus du point $s\in S$ et si ${\cal{F}}$ est un faisceau
sur $X$, on note  ${\cal{F}}_{s}$ la restriction de ${\cal{F}}$ \`a $X_{s}$.

Si $X$ est une vari\'et\'e de Cohen-Macaulay on d\'esigne par
$\omega_{_{X}}$ le faisceau dualisant. Si ${\cal{F}}$ est un faisceau de
Cohen-Macaulay de dimension $d$ sur une vari\'et\'e de Cohen-Macaulay de
dimension $n$, on note ${\cal{F}}^{\vee}$ le {\em $\omega$-dual}
de ${\cal{F}}$ \ie le faisceau
$\ul{\Ext}^{n-d}_{{\cal{O}}_{X}}({\cal{F}},\omega_{_{X}})$.

Si $(L^{\cdot},d_{L})$ est un complexe de
${\cal{O}}_{S}$-modules, on d\'esigne par $L^{\cdot}[m]$  le complexe
translat\'e de $m$ places \`a gauche, de diff\'erentielle $(-1)^{m}d_{L}$ et
par $\tau_{\leq m}(L^{\cdot})$ le sous-complexe de $L$ d\'efini par
$\dots\lra L^{m-2}\lra L^{m-1}\lra \Ker(d^{m})\lra 0$.

Si $(K^{\cdot},d_{K})$ est un autre complexe de
${\cal{O}}_{S}$-modules, le complexe de ${\cal{O}}_{S}$-modules
$\ul{\Hom}^{\cdot}(L^{\cdot},K^{\cdot})$ est d\'efini en degr\'e $n$ par
$\displaystyle \prod_{-p+q=n}\ul{\Hom}(L^{p},K^{q})$, de diff\'erentielle
$d(f)^{p}=d_{K}\circ f^{p}-(-1)^{n}f^{p+1}d_{L}$. Finalement, on note
$L^{\cdot*}$ le  complexe $\ul{\Hom}(L^{\cdot},{\cal{O}}_{S})$.

Par $D(S)$ on d\'esigne la cat\'egorie d\'eriv\'ee de la cat\'egorie des
${\cal{O}}_{S}$-modules, par
$D_{c}(S)$ la sous-cat\'egorie pleine des complexes
\`a cohomologie coh\'erente, par $D^{b}(S)$ la sous-cat\'egorie pleine des
complexes born\'es \`a gauche et \`a droite et enfin, par $D^{b}_{c}(S)$
l'intersection dans $D(S)$ de $D_{c}(S)$ avec $D^{b}(S)$.

\bigskip\medskip
{\em Remerciements:} Je tiens \`a remercier Joseph Le Potier pour son
aide pendant la pr\'eparation du pr\'esent travail.

\np
\tableofcontents
\np

\section{Le module $\Theta_{\planp}(d)$ des th\^eta-caract\'eristiques
des courbes planes}

Nous allons rappeler les d\'efinitions
et r\'esultats de \cite{11} sur l'espace de modules des
th\^eta-caract\'eristiques semi-stables et celles
de \cite{9} sur le module des faisceaux semi-stables de dimension 1 sur le plan
projectif dont on a besoin dans la suite.

\subsection{Th\^eta-caract\'eristiques des courbes planes }

Soit ${\cal{F}}$ un faisceau pur de dimension $1$ sur le plan
projectif,
\ie tel que tout sous-faisceau non nul de ${\cal{F}}$ soit de
dimension $1$. C'est un faisceau de Cohen-Macaulay; il admet donc
une r\'esolution  $$0\lra A\hfl{\alpha}{}B\lra{\cal{F}}\lra 0$$ par
des faisceaux localement libres sur $\planp$. On appelle {\em
support sch\'ematique} et on note
$\supps({\cal{F}})$ la courbe d\'efinie par $\det(\alpha)$. C'est une
courbe de Gorenstein dont le degr\'e est \'egale \`a la {\em
multiplicit\'e} de
${\cal{F}}$, not\'ee $\mult({\cal{F}})$, \ie au coefficient directeur
du polyn\^ome de Hilbert de
${\cal{F}}$.   Posons
${\cal{F}}^{\vee}=\ul{\Ext}_{{\cal{O}}_{\planp}}^{1}
          ({\cal{F}},\omega_{_{\planp}}).$
          Si l'on consid\`ere
${\cal{F}}$ comme faisceau de modules sur son support sch\'ematique,
on a ${\cal{F}}^{\vee}\simeq
 {\cal{F}}^{*}\otimes_{{\cal O}_{C}}\omega_{_{C}}$, o\`u
${\cal{F}}^{\vee}=
 \ul{\Hom}_{{\cal{O}}_{C}}({\cal{F}},\omega_{_{C}})$.
De plus, on a $\chi({\cal{F}}^{\vee})=-\chi({\cal{F}})$ et le
morphisme canonique d'\'evaluation
${\cal{F}}\lra{\cal{F}}^{\vee\vee}$ est un isomorphisme (\cite{11}, 1.1).

\begin{defi}\label{theta-def} On appelle th\^eta-caract\'eristique sur
$\planp$ la donn\'ee d'un couple $({\cal{F}},\sigma)$ form\'e d'un
${\cal{O}}_{\planp}$-module coh\'erent de dimension
$1$ et d'un isomorphisme sym\'etrique
$\sigma:{\cal{F}}\lra{\cal{F}}^{\vee}.$
\end{defi}

Si $({\cal{F}},\sigma)$ est une th\^eta-caract\'eristique, alors
$\chi({\cal{F}})=0$. Un {\em morphisme} de th\^eta-caract\'eristiques
$\varphi:({\cal{E}},\tau)\lra({\cal{F}},\sigma)$ est la donn\'ee d'un
morphisme de ${\cal{O}}_{\planp}$-modules
$\varphi:{\cal{E}}\lra{\cal{F}}$ tel que
$\tau=^{\vee}\varphi\circ\sigma\circ\varphi$. Soit
${\cal{E}}\subset{\cal{F}}$ un sous-faisceau coh\'erent de
${\cal{F}}$. On appelle {\em $\sigma$-orthogonal} de ${\cal{E}}$ et
on note ${\cal{E}}^{\ort}$ le noyau du morphisme compos\'e
${\cal{F}}\lra{\cal{F}}^{\vee}\lra{\cal{E}}^{\vee}.$ Un
sous-faisceau coh\'erent ${\cal{E}}\subset{\cal{F}}$ est dit  {\em
$\sigma$-isotrope} si ${\cal{E}}\cap{\cal{E}}^{\ort}\not=(0)$. Si
l'on a ${\cal{E}}\subset{\cal{E}}^{\ort}$ on dit que ${\cal{E}}$
est  {\em totalement $\sigma$-isotrope}.

Une th\^eta-caract\'eristique
$({\cal{F}},\sigma)$ est dite {\em semi-stable (resp. stable)} si
pour tout sous-faisceau totalement $\sigma$-isotrope non nul
${\cal{E}}\subset{\cal{F}}$ on a
$$\chi({\cal{E}})\le 0\text{(resp. $<$)}.$$ Une th\^eta-caract\'eristique
$({\cal{F}},\sigma)$ est semi-stable si et seulement si son faisceau
sous-jacent ${\cal{F}}$ est semi-stable. Une th\^eta-caract\'eristique stable
est somme directe orthogonale de th\^eta-caract\'eristiques dont les faisceaux
sous-jacents sont stables et 2 \`a 2 non-isomorphes  (\cite{11},
Prop. 1.4). On dit qu'une th\^eta-caract\'eristique est ${\cal{O}}$-{\em
stable} si son faisceau sous-jacent est stable.

La cat\'egorie des th\^eta-caract\'eristiques semi-stables
n'est pas ab\'elienne.
On
peut n\'eanmoins d\'efinir une notion de
${\cal{S}}$-\'equivalence de deux faisceaux quadratiques: Soit
$({\cal{F}},\sigma)$ une th\^eta-caract\'eristique semi-stable. Alors
il existe une filtration par des sous-faisceaux coh\'erents totalement
isotropes
$ 0\subset{\cal{F}}_{1}
 \subset\dots\subset{\cal{F}}_{\ell}\subset{\cal{F}}
$ telle que les faisceaux
${\cal{E}}_{i}={\cal{F}}_{i}/{\cal{F}}_{i-1}$, avec
${\cal{F}}_{0}=0$, soient stables de caract\'eristique
d'Euler-Poincar\'e nulle pour $i=1,\dots,\ell$ et telle que
$({\cal{F}}_{\ell}^{\ort}/{\cal{F}}_{\ell},\sigma)$ soit une
th\^eta-caract\'eristique stable.  Une th\^eta-caract\'eristique
$({\cal{H}},\sigma)$ est dite {\em hyperbolique} si elle est
isomorphe \`a une th\^eta-caract\'eristique
$({\cal{G}}\osum{\cal{G}}^{\vee},\tau)$, o\`u ${\cal{G}}$ est stable
en tant que faisceau et o\`u $\tau$ est donn\'e par la matrice
suivante:
$$
\bordermatrix{%
                & {\cal{G}} & {\cal{G}}^{\vee} \cr
{\cal{G}}^{\vee}&     0     &        1  \cr
       {\cal{G}}&     1     &        0  \cr }
$$ Munissons ${\cal{H}}_{i}={\cal{E}}_{i}\osum{\cal{E}}_{i}^{\vee}$
de la structure hyperbolique $\sigma_{i}=\tau$. On d\'efinit la
th\^eta-caract\'eristique gradu\'e associ\'ee \`a  $({\cal{F}},\sigma)$ comme
la somme directe orthogonale suivante:
$$
\gr({\cal{F}},\sigma)=
 \left({\cal{F}}_{\ell}^{\ort}/{\cal{F}}_{\ell},\sigma\right)\osum
 \left(\osum_{i=1}^{\ell}({\cal{H}}_{i},\sigma_{i})\right).
$$ La filtration n'est pas unique. Le gradu\'e par contre, ne d\'epend
pas, \`a isomorphisme pr\`es, de la filtration choisie. Par cons\'equent,
on dira que deux th\^eta-caract\'eristiques
$({\cal{E}},\tau)$ et $({\cal{F}},\sigma)$ sont {\em
${\cal{S}}$-\'equivalentes} si leurs gradu\'es associ\'es respectifs sont
isomorphes. On dira que $({\cal{F}},\sigma)$ est {\em poly-stable} si
$$({\cal{F}},\sigma)=\left(\osum_{i}{\cal{H}}_{i}^{m_{i}}\right)
\osum\left(\osum_{j}{\cal{F}}_{j}^{n_{j}}\right),$$
la somme directe \'etant orthogonale, les ${\cal{H}}_{i}$ (et
de m\^eme les ${\cal{F}}_{j}$) \'etant deux \`a deux non-isomorphes et
${\cal{H}}_{i}={\cal{E}}_{i}\osum{\cal{E}}_{i}^{\vee}$ \'etant hyperbolique
{\em stricte} \ie ${\cal{E}}_{i}\not\simeq{\cal{E}}_{i}^{\vee}$.
D'apr\`es ce qui pr\'ec\`ede le gradu\'e d'une th\^eta-caract\'eristique est
poly-stable
et une th\^eta-caract\'eristique poly-stable est isomorphe \`a son gradu\'e.
Les points ferm\'es du module des th\^eta-caract\'eristiques qui sera d\'efini
dans
la section suivante seront exactement les th\^eta-caract\'eristiques
poly-stables.

\subsection{La construction du module des th\^eta-caract\'eristiques sur
$\planp$.}\label{thetaconstruction}

Soit $S$ une vari\'et\'e alg\'ebrique, ${\cal{F}}$ un
${\cal{O}}_{S\times \planp}$-module, $S$-plat. On pose
$${\cal{F}}^{\vee}=\ul{\Ext}_{{\cal{O}}_{S\times \planp}}^{1}
        ({\cal{F}},pr_{2}^{*}\omega_{_{\planp}}).$$ Si ${\cal{F}}$
est une famille $S$-plate de faisceaux purs de dimension $1$, alors
${\cal{F}}^{\vee}$ est $S$-plat, le morphisme canonique
${\cal{F}}\lra{\cal{F}}^{\vee\vee}$ est un isomorphisme et l'on a
$({\cal{F}}^{\vee})_{s}\simeq({\cal{F}}_{s})^{\vee}$ (\cite{11}, 3.1). Une {\em
famille de th\^eta-caract\'eristiques}
para\-m\'etr\'ee par la vari\'et\'e
alg\'ebrique $S$ est la donn\'ee d'un couple
$({\cal{F}},\sigma)$, form\'e d'un ${\cal{O}}_{S\times X}$-module
coh\'erent
${\cal{F}}$, $S$-plat, et d'un isomorphisme sym\'etrique
$\sigma:{\cal{F}}\lra{\cal{F}}^{\vee}$.

Consid\'erons le foncteur
$\underline{\Theta}_{\planp}(d)$ associant \`a la vari\'et\'e alg\'ebrique
$S$ l'ensemble des classes d'iso\-morphisme de familles de
th\^eta-caract\'eristiques $({\cal{F}},\sigma)$, param\'etr\'ees par $S$,
telles que, pour tout point ferm\'e $s$ de $S$, la
th\^eta-caract\'eristique $({\cal{F}}_{s},\sigma_{s})$ soit semi-stable
et de polyn\^ome de Hilbert $P(m)=md$. Pour le foncteur
$\underline{\Theta}_{\planp}(d)$ il existe un espace de modules
grossier, not\'e $\Theta_{\planp}(d)$. C'est une vari\'et\'e projective
dont l'ensemble des points ferm\'es est l'ensemble des classes de
${\cal S}$- \'equivalence de th\^eta-caract\'eristiques semi-stables de
degr\'e
$d$ sur $\planp$. De plus, il existe un ouvert non-vide
$\Theta_{\planp}^{s}(d)\subset\Theta_{\planp}(d)$ dont les points
repr\'esentent les classes d'isomorphisme des  th\^eta-caract\'eristiques
stables (\cite{11}, 0.2).

Dans ce qui suit on a besoin d'expliciter la construction de
$\Theta_{\planp}(d)$:
soit $N$ un entier tel que pour toute th\^eta-caract\'eristique semi-stable
$({\cal{F}},\sigma)$ le faisceau ${\cal{F}}(N)$ soit engendr\'e par
ses sections globales et de cohomologie sup\'erieure nulle.
Ceci est possible, la famille des th\^eta-caract\'eristiques semi-stables de
degr\'e
fix\'e \'etant limit\'ee (\cite{11}, Prop. 3.3]). Soit
$H=k^{dN}$ et
${\tr{H}}=H\otimes_{k}{\cal{O}}_{\planp}(-N)$.
Consid\'erons, pour toute vari\'et\'e alg\'ebrique
$S$, les triplets $({\cal{F}},\alpha,\sigma)$ form\'es d'un quotient
coh\'erent
${\tr{H}}_{S}\hfl{\alpha}{}{\cal{F}}$, $S$-plat, et d'un
isomorphisme sym\'etrique $\sigma:{\cal{F}}\lra{\cal{F}}^{\vee}$
satisfaisant aux conditions suivantes: pour tout point ferm\'e $s$ de
$S$, $({\cal{F}}_{s},\sigma_{s})$ est une th\^eta-caract\'eristique
semi-stable de multiplicit\'e $d$ et $\alpha$ induit un isomorphisme
$H\otimes{\cal{O}}_{S}\simeq pr_{1*}({\cal{F}}(N))$. Deux triplets
$({\cal{F}},\alpha,\sigma)$ et
$({\cal{F}}',\alpha',\sigma')$ sont dits \'equivalents s'il existe un
isomorphisme $\varphi:{\cal{F}}\lra{\cal{F}}'$ tel que
$\alpha'=\varphi\circ \alpha$ et tel que
$\sigma'\varphi=(^{\vee}\varphi)^{-1}\circ\sigma$. Notons par
$[{\cal{F}},\alpha,\sigma]$  la classe d'\'equivalence du triplet
$({\cal{F}},\alpha,\sigma)$ et par
$\underline{T}^{ss}(d,N)(S)$ l'ensemble des classes d'\'equivalence
de tels triplets. On obtient un foncteur
$$
\underline{T}^{ss}(d,N):\text{Vari\'et\'es alg\'ebriques}
\lra\text{Ensembles},
$$ dont on montre la repr\'esentabilit\'e (\cite{11}, 7.3).
Notons  $T^{ss}(d,N)$ le sch\'ema qui le repr\'esente. Le groupe
$GL(H)$ op\`ere sur $T^{ss}(d,N)$ en associant \`a $g\in GL(H)$ et au
triplet
$[{\cal{F}},\alpha,\sigma]$ le triplet $[{\cal{F}},\alpha\circ
g^{-1},\sigma]$. D'apr\`es \cite{11}, l'espace de modules
$\Theta_{\planp}(d)$ s'identifie au quotient de Mumford de
$T^{ss}(d,N)/GL(H)$.

Le stabilisateur d'un point $[{\cal{F}},\alpha,\sigma]$ s'identifie
au groupe des automorphismes de la th\^eta-caract\'eristiques
$({\cal{F}},\sigma)$. En effet, soit
$\chi\in\Aut({\cal{F}},\sigma)$. De $\chi$ on d\'eduit un
isomorphisme $\H^{0}({\cal{F}}(N))\lra\H^{0}({\cal{F}}(N))$
induisant un isomorphisme $\overline\chi\in GL(H)$:
$$
\begin{diagram}
\H&\hfl{\alpha}{}&\H^{0}({\cal{F}}(N))\\
\vfl{\overline\chi}{}&&\vfl{\chi}{}\\
\H&\hfl{\alpha}{}&\H^{0}({\cal{F}}(N))
\end{diagram}
$$ qui est dans le stabilisateur. R\'eciproquement, un \'el\'ement
$g\in GL(H)$ du stabilisateur induit, comme
$[{\cal{F}},\alpha,\sigma] =[{\cal{F}},\alpha\circ g^{-1},\sigma]$,
un isomorphisme $\phi:{\cal{F}}\lra{\cal{F}}$ respectant $\sigma$.
Ces deux op\'erations sont inverses l'une de l'autre, d'o\`u
l'identification voulue.

La vari\'et\'e $T^{ss}(d,N)$ est une vari\'et\'e lisse (\cite{11},
Prop. 7.4) et
$GL(H)/\{\pm Id\}$ op\`ere librement sur l'ouvert des
th\^eta-caract\'eristiques
${\cal{O}}$-stables. Il en d\'ecoule que $\Theta_{\planp}(d)$ est
normale et n'a que des singularit\'es rationnelles. Il en d\'ecoule
aussi que l'ouvert
$\Theta_{\planp}^{os}(d)$ des th\^eta-caract\'eristiques
${\cal{O}}$-stables est lisse. La codimension du ferm\'e des
th\^eta-caract\'eristiques semi-stables, non ${\cal{O}}$-stables est au
moins
$2$ (\cite{11}, Prop. 8.6). On utilisera ce r\'esultat plus
tard sans le rappeler explicitement.

\subsection{La structure sch\'ematique de $\Theta^{r}(d)$.}
\label{BNlocus}

Soit $({\cal{F}},\sigma)$ une famille de th\^eta-caract\'eristiques
param\'etr\'ee
par la vari\'et\'e alg\'ebrique $S$. D\'efinissons les ensembles
$$S^{r}=\{s\in S/\h^{0}({\cal{F}}_{s})\geq r+1,
\h^{0}({\cal{F}}_{s})= r+1\bmod 2\}.$$
D'apr\`es le th\'eor\`eme de semi-continuit\'e et le th\'eor\`eme d'invariance
mod
$2$ les $S^{r}$ sont ferm\'es.
Pour la d\'efinition de la structure sch\'ematique des $S^{r}$, rappelons qu'on
appelle approximation  de le cohomologie de ${\cal{F}}$
un complexe
$$E^{\cdot}\text{: }0\lra E^{0}\efl{d}{} E^{1}\lra 0$$, form\'e de
${\cal{O}}_{S}$-modules coh\'erents
localement libres, tel que
pour tout changement de base $S'\hfl{g}{}S$,
$$
\begin{diagram}
S'\times\planp&\hfl{g'}{}&S\times\planp\\
\vfl{f'}{}&&\vfl{}{f}\\
S'&\hfl{g}{}&S
\end{diagram}
$$
on ait $\H^{i}(g^{*}E^{\cdot})=\R^{i}f'_{*}g^{\prime *}{\cal{F}}$. En
particulier, pour $s\in S$, on a
$\H^{i}(E^{\cdot}_{s})=\H^{i}(X_{s},{\cal{F}}_{s})$.
Une telle approximation est dite anti-sym\'etrique si
$E^{1}=E^{0*}$ et si $d$ est anti-sym\'etrique.
Une approximation de la cohomologie globale de ${\cal{F}}$ existe
toujours; par contre pour une
approximation anti-sym\'etrique on n'a que l'existence, pour tout point
$s\in S$, d'un voisinage ouvert $U_{s}$ de $s$ telle que sur $U_{s}$ il
existe une approximation anti-sym\'etrique
$$0\lra Z^{0}\efl{d}{}Z^{*}\lra 0$$
de ${\cal{F}}$ avec $Z$ de rang $z=\h^{0}({\cal{F}}_{s})$ et $d(s)=0$.
Localement sur $U$, la structure sch\'ematique de $S^{r}$ est
d\'efinie par les mineurs d'ordre $(z-r)$ de $\alpha$. L'id\'eal engendr\'e par
ces mineurs ne d\'epend pas de l'approximation anti-sym\'etrique particuli\`ere
choisie. De ce fait, cette construction
d\'efinit donc un faisceau d'id\'eaux sur $S$ entier d\'efinissant ainsi la
structure sch\'ematique de $S^{r}$.

En particulier, si $s\in S^{r}\backslash S^{r+1}$, le sous-sch\'ema $S^{r}$ est
 d\'efini au voisinage de $s$ par des $1$-mineurs d'une matrice
anti-sym\'etrique,
 donc par $\frac{r(r-1)}{2}$ \'equations. Ceci d\'emontre, en utilisant le
lemme de
Krull, que
 $$\codim_{S}(S^{r})\leq\frac{r(r-1)}{2},$$
si $S^{r}\backslash S^{r+1}$ est non-vide.

Consid\'erons maintenant la vari\'et\'e $T^{ss}(d)$ de la section
\ref{thetaconstruction} et le bon quotient
$\varphi:T^{ss}(d)\lra\Theta(d)$. Sur $T^{ss}(d)$ on a une famille
universelle $({\cal{F}},\sigma,\alpha)$. Soient les sous-sch\'emas
$T^{ss,r}(d)\subset T(d)$ d\'efinis comme ci-dessus. Ces sous-sch\'emas sont
invariants
sous l'action de $GL(H)$. On d\'efinit
$\Theta^{r}(d)$ comme \'etant l'image de $T^{ss,r}(d)$. D'apr\`es les
propri\'et\'es
de bons quotients, $\Theta^{r}(d)$ est un sous-sch\'ema ferm\'e de $\Theta(d)$.
L'ensemble sous-jacent \`a $\Theta^{r}(d)$ s'identifie \`a
$$\Theta^{r}(d)=\{[{\cal{F}},\sigma]\in\Theta(d)/\h^{0}({\cal{F}})\ge r+1,
\h^{0}({\cal{F}})=(r+1)\bmod 2\}.$$

\subsection{L'espace de modules $N_{\planp}(d,0)$ de Le
Potier}\label{LePotierspace}

Soit $\ul{N}_{\planp}(d,\chi)$ le foncteur qui associe \`a la vari\'et\'e
alg\'ebrique $S$ l'ensemble des classes d'isomorphismes de familles
$S$-plates de faisceaux semi-stables de dimension $1$ de degr\'e
$d$ et de caract\'eristique d'Euler-Poincar\'e $\chi$. D'apr\`es
\cite{10}, il existe un espace de modules grossier projectif, not\'e
$N_{\planp}(d,\chi)$, pour ce foncteur
et l'ensemble de ses points ferm\'es s'identifie aux classes de
${\cal{S}}$-\'equivalence de faisceaux semi-stables. De plus, on a un
morphisme, appel\'e morphisme support sch\'ematique,
$$\sigma_{N}:N(d,\chi)_{\planp}\lra\proj^{M},$$ avec $M=d(d+3)/2$, qui
associe \`a un faisceau semi-stable de dimension $1$ son support
sch\'ematique.

D'apr\`es \cite{9},
la vari\'et\'e
$N_{\planp}(d,\chi)$ est irr\'eductible, normale et localement
factorielle. De plus, son groupe de Picard est isomorphe \`a
\`a un groupe ab\'elien \`a deux g\'en\'erateurs. Pour la description de ces
g\'en\'erateurs dont
on a besoin dans ce qui suit, on va rappeler la construction de cet
espace de modules.
Elle se fait de la mani\`ere suivante: d'abord on montre
que la famille des faisceaux semi-stables de polyn\^ome de Hilbert
fix\'e est limit\'ee. Ensuite, on choisit un entier $N$ de fa\c con que
pour  tout faisceau semi-stable ${\cal{F}}$ de dimension $1$, de
multiplicit\'e
$d$  et de caract\'eristique d'Euler-Poincar\'e $\chi$, le faisceau
${\cal{F}}(N)$  soit engendr\'e par ses sections, et que
$\H^1({\cal{F}}(N))=0.$  Un tel entier \'etant choisi, on consid\`ere
un espace vectoriel $\H$  de dimension
$n=dN+\chi$ et le fibr\'e vectoriel ${\tr{H}}= H_{k}\otimes{\cal
O}_{\planp}(-N)$.  Le groupe $SL(H)$ op\`ere de mani\`ere naturelle sur
le sch\'ema de Hilbert-Grothendieck
$\Groth^{d,\chi}({\tr{H}})$  des faisceaux quotients de
${\tr{H}}$  de multiplicit\'e
$d$
 et de caract\'eristique d'Euler-Poincar\'e $\chi$. Quitte \`a choisir un
entier $m$ suffisamment grand, ce sch\'ema se plonge dans une
grassmannienne par le plongement de Grothendieck: au quotient
${\tr{H}}\lra{\cal{F}}$ on associe le quotient
$H\otimes_{k}\H^{0}({\cal{O}}_{\planp}(m-N))\lra{\cal{F}}(m)$. De
plus, l'action de $SL(H)$ se rel\`eve en une action lin\'eaire sur la
grassmannienne. Maintenant, quitte \`a choisir $N$ et ensuite $m$
suffisamment grand, l'ouvert $\Omega^{ss}$ des points semi-stables pour
l'action de $SL(H)$ correspond aux faisceaux quotients de
${\tr{H}}$ qui sont semi-stables et tels que le morphisme
d'\'evaluation
$\H\ra\H^0({\cal{F}}(N))$ soit un isomorphisme \cite{10}. Alors
$N_{\planp}(d,\chi)$ s'identifie au quotient de Mumford
$$N_{\planp}(d,\chi)=\Omega^{ss}/SL(\H).$$

Venons-en \`a la description des g\'en\'erateurs de
$\Pic(N_{\planp}(d,0))$: ce sont les fibr\'e inversible
${\cal{L}}_{N}$, image r\'eciproque de
${\cal{O}}_{\proj^{M}}(1)$ sous le morphisme
$\sigma_{N}$, et le fibr\'e inversible ${\cal{D}}_{N}$, appel\'e fibr\'e
d\'eterminant, associ\'e au diviseur ``th\^eta'' des faisceaux
semi-stables ${\cal{F}}$ ayant au moins une section (\cite{9}, 1.1). On a une
autre
description, davantage fonctorielle, de
ces g\'en\'erateurs:  consid\'erons pour cela l'alg\`ebre  de Grothendieck
$K(\planp)$ de $\planp$. Elle est muni d'une forme quadratique \`a
valeurs enti\`ere, d\'efinie en associant \`a $u\in K(\planp)$ l'entier
relatif $\chi(u^{2})$.  Cette forme provient d'une forme
quadratique sur
$K_{top}(\planp)$. On note $<\ ,\ >$ la forme bilin\'eaire associ\'e.
Si $x$ (resp. $y$) est de rang $r$, de degr\'e $d$ et de
caract\'eristique d'Euler-Poincar\'e $\chi$ (resp. $r',d',\chi'$), on
a:
$$<x,y>=r\chi'+r'\chi+dd'-rr'.$$ Soit ${\cal{F}}$ une famille
$S$-plate de faisceaux semi-stables de dimension 1 sur $\planp$,
param\'etr\'ee par la vari\'et\'e alg\'ebrique $S.$  Supposons de plus que le
groupe alg\'ebrique $G$ op\`ere sur $S$ et ${\cal{F}}$.
 Consid\'erons le diagramme
$$\begin{diagram}
S\times X&\hfl{p}{}&X\\
\sfl{q}{}&&\cr S&&\\
\end{diagram}
$$  o\`u $p$  et $q$ sont les projections
canoniques. Alors \`a $u\in K(\planp)$, on peut associer l'\'el\'ement
$${\cal{L}}_{{\cal{F}}}(u)=\det(q_{!}({\cal F}.p^*(u))$$ du groupe
$\Pic^{G}(S)$ des fibr\'e inversibles sur $S$ muni d'une action de
$G$. On peut appliquer cette construction en particulier au
faisceau quotient universel sur $\Omega^{ss}\times X$. Le fibr\'e
inversible ${\cal{L}}_{{\cal{F}}}(u)$ est alors muni d'une action
de $GL(H)$. D\'esignons par $Z(d,\chi)$ l'orthogonal de
$(d,\chi)$ dans $K(\planp)$ par rapport \`a la forme quadratique
d\'efinie ci-dessus. Le r\'esultat de \cite{8} affirme
maintenant que pour $u\in  Z(d,\chi)$  il existe dans
$\Pic(N_{\planp}(d,\chi))$ un \'el\'ement ${\cal{L}}(u)$  caract\'eris\'e
par la propri\'et\'e universelle suivante: pour toute famille plate
${\cal{F}}$  de faisceaux semi-stables sur $X$  de degr\'e $d$ et de
caract\'eristique d'Euler-Poincar\'e $\chi$ param\'etr\'ee par la vari\'et\'e
alg\'ebrique $S$ on a $${\cal{L}}_{{\cal F}}(u)
=f_{{\cal{F}}}^*({\cal{L}}(u))$$  o\`u $f_{{\cal F}}:S\lra
N_{\planp}(d,\chi)$ est le morphisme modulaire associ\'e \`a la famille
${\cal{F}}$.

Pour les g\'en\'erateurs ${\cal{L}}_{N}$ et ${\cal{D}}_{N}$ on a alors
${\cal{L}}_{N}={\cal{L}}(u)$ avec $u$ la classe d'un point,
et
${\cal{D}}_{N}={\cal{L}}(u)$ avec $u=-[{\cal{O}}_{\planp}]$.

\np
\section{L'ouvert des th\^eta-caract\'eristiques ne
passant pas par un point.}\label{Ua}

Soit $U_{a}$ l'ouvert correspondant aux th\^eta-caract\'eristiques
semi-stables ne passant pas par le point $a\in\planp$.
On note $U_{a,p}(d)$ (resp. $U_{a,i}(d)$) le ferm\'e de
$U_{a}(d)$ correspondant aux th\^eta-caract\'eristiques paires (resp.
impaires).
Soit de plus
$$U^{r}_{a}(d):=
\{[{\cal{F}},\sigma]\in U_{a}/\h^{0}({\cal{F}})\ge r+1,
\h^{0}({\cal{F}})=(r+1)\bmod 2\},$$
qu'on muni de sa structure de vari\'et\'e d\'eterminantielle naturelle.
 L'objet de cette
section est de d\'emontrer le th\'eor\`eme suivant:

\begin{th}\label{codimU} L'ouvert $U_{a}(d)$ est, si $d$ est pair, r\'eunion de
deux composantes irr\'eductibles correspondant aux th\^eta-caract\'eristiques
paires et impaires; si $d$ est impair, il est r\'eunion de trois composantes
irr\'eductibles correspondant aux
th\^eta-caract\'eristiques paires, impaires et
canoniques. De plus, le ferm\'e $U^{1}_{a}(d)$ est irr\'eductible de
codimension $1$
pour $d\ge 5$, vide sinon. Le ferm\'e $U^{2}_{a}(d)$
est de codimension au moins deux sauf si $d=6$ o\`u il est
irr\'eductible de codimension $1$.
\end{th}

Ceci donnera,
en corollaire, apr\`es l'\'etude du compl\'ementaire $F_{a}(d)$ de
$U_{a}(d)$ le th\'eor\`eme \ref{BN} et aussi une nouvelle preuve de
l'irr\'eductibilit\'e des composantes paires et impaires
de \cite{11}.

Soit $\proj_{1}$ une droite du plan projectif ne passant
pas par ce point. Consid\'erons la projection de centre a sur la droite:
$$\pi:\planp-\{a\}\lra\proj_{1}$$
Les images directes sous $\pi$ des th\^eta-caract\'eristiques ne
passant pas par
$a$ se d\'ecrivent en termes de fibr\'es
$\omega_{_{\proj_{1}}}$-quadratiques sur la droite projective muni d'un
morphisme $\varphi:G\lra G(1)$ compatible avec la structure quadratique:

\subsection{Fibr\'es de Higgs quadratiques}
Rappelons bri\`evement quelques g\'en\'eralit\'es sur les fibr\'es de Higgs:
soit $C$ une courbe projective lisse et soit $L$ un fibr\'e inversible sur
$C$. On appelle {\em faisceau de Higgs} la donn\'ee d'un
couple $(G,\varphi)$ form\'e d'un faisceau coh\'erent $G$ et d'un
morphisme $\varphi:G\lra G\otimes L$. Si $G$ est localement libre on dira
{\em fibr\'e de Higgs}. Un morphisme de faisceaux de Higgs de
$F'=(G',\varphi')$ dans $F''=(G'',\varphi'')$ est la donn\'ee d'un morphisme
de faisceau $f:G'\lra G''$ tel que le diagramme suivant commute:
$$\begin{diagram} G'&\efl{f}{}&G''\\
\sfl{\varphi'}{}&&\sfl{\varphi''}{}\\
G'\otimes L&\efl{f\otimes id}{}&G''\otimes L\\
\end{diagram}
$$
Ainsi un sous-faisceau de Higgs est la donn\'ee d'un sous-faisceau $G'\subset
G$ tel que $\varphi(G')\subset G'\otimes L$.
Un fibr\'e de Higgs
est dit {\em semi-stable} (resp. {\em stable}) si pour tout
sous-faisceau propre de Higgs non nul $G'$ de
$F=(G,\varphi)$ on a
$$\mu(G')\le\mu(G) \text{ (resp. $<$)}.$$
Comme dans le cas des fibr\'es vectoriels il suffit de v\'erifier
cette condition pour les {\em sous-fibr\'es} de Higgs. Si
$F=(G,\varphi)$ est un faisceau de Higgs on note  $H^{i}(C,F)$ la
cohomologie de
$F$,
\ie l'hypercohomologie du complexe
$$0\lra G\efl{\varphi}{}G\otimes L\lra 0.$$ Si
$F'=(G',\varphi')$ et $F''=(G'',\varphi'')$ sont deux faisceaux de Higgs le
faisceaux des homomorphismes de $F'$ dans $F''$ est d\'efini par
$(\ul{\Hom}(G',G''),\varphi)$ avec
$\varphi(s)=\varphi''s-s\varphi'$ o\`u $s$ est une section locale de
$\ul{\Hom}(G',G'')$. On note $\Ext^{q}(F',F'')$ la cohomologie
de $\ul{\Hom}(F',F'')$. Bien s\^ur, on a
$\Ext^{0}(F',F'')=\Hom(F',F'')$. Le $\omega_{_{C}}$-dual d'un faisceau de
Higgs $F=(G,\varphi)$ est d\'efini par la paire
$(G^{\vee},\psi)$ o\`u $\psi:G^{\vee}\lra G^{\vee}\otimes L$ est d\'efini par
la transpos\'ee $^{\vee}\varphi:G^{\vee}\otimes L^{*}\lra G^{\vee}$ de
$\varphi$. On le note $F^{\vee}$. Si $F$ est un fibr\'e de Higgs, le morphisme
canonique d'\'evaluation est un isomorphisme.

\begin{defi}
On appelle fibr\'e de Higgs quadratique la donn\'ee d'un fibr\'e de Higgs $F$
et d'un isomorphisme de Higgs sym\'etrique $\sigma:F\lra F^{\vee}$.
\end{defi}
Soit $(F,\sigma)$ un fibr\'e de Higgs quadratique, $F'\subset F$ un
sous-faisceau (de Higgs). On d\'efinit {\em l'orthogonal} de $F'$ et
l'on note
$F^{\prime\ort}$ le noyau du morphisme compos\'e $$F\lra F^{\vee}\lra
F^{\prime\vee}.$$ Un sous-faisceau $F'$ est dit isotrope si $F'\cap
F^{\prime\ort}\not=0$ et totalement isotrope si $F'\subset F^{\prime\ort}$.
Un fibr\'e de Higgs quadratique est dit {\em semi-stable} (resp. {\em
stable}) si pour tout sous-faisceau totalement isotrope
$F'=(G',\varphi)$ de $F=(G,\varphi)$ on a
$$\mu(G')\le\mu(G) \text{ (resp. $<$}).$$
\begin{lemme}\label{ss-ss}
Un fibr\'e de Higgs quadratique est semi-stable si et seulement si
son fibr\'e de Higgs sous-jacent est semi-stable.
\end{lemme}
\begin{proof}
Analogue \`a la proposition 1.4 de \cite{11}.
\end{proof}
Cette \'equivalence n'est pas vraie pour la notion de stabilit\'e. On
dira qu'un fibr\'e de Higgs quadratique est {\em ${\cal{O}}$-stable}
si son fibr\'e de Higgs sous-jacent est stable.

Soit $(G,\sigma)$ un fibr\'e vectoriel $\omega_{_{C}}$-quadratique.
L'isomorphisme $\sigma$ d\'efinit une involution $\imath$ sur
l'espace vectoriel $\Ext^{1}(G,G^{\vee})$. Le sous-espace vectoriel
des \'el\'ements sym\'etriques (resp. antisym\'etriques) est not\'e
$\Ext^{1}_{sym}(G,G^{\vee})$ (resp. $\Ext^{1}_{asym}(G,G^{\vee})$).
En notant $\Ext^{1}_{sym}(G,G)$ (resp.
$\Ext^{1}_{asym}(G,G)$) on suppose que l'on a identifi\'e $G$ et
$G^{\vee}$ via $\sigma$. \begin{defi}
Un fibr\'e vectoriel $\omega_{_{C}}$-quadratique $(G,\sigma)$ est dit rigide si
$$\Ext^{1}_{asym}(G,G)=0.$$
\end{defi}
Un fibr\'e de Higgs quadratique dont le fibr\'e vectoriel quadratique
sous-jacent est rigide est dit {\em quasi-rigide}.

\subsubsection{Le cas de la droite projective.}

Soit maintenant $C=\proj_{1}$ et $L={\cal{O}}_{\droitep}(1)$.
\'Etant donn\'e un
fibr\'e vectoriel $G$ sur $\droitep$ on \'ecrit
$G=\osum_{\ell\in\reln}r_{\ell}{\cal{O}}_{\proj_{1}}(\ell)$. La
suite des entiers $r_{\ell}\not=0$ s'appelle {\em spectre} de $G$,
les entiers $r_{\ell}$ s'appellent multiplicit\'es.
\begin{lemme} Soit $(G,\varphi,\sigma)$ un fibr\'e de Higgs quadratique
semi-stable. Alors le spectre de $G$ est connexe.
\end{lemme}
\begin{proof} En raison du lemme \ref{ss-ss} on peut supposer que
$(G,\varphi)$ est semi-stable, puis on applique le lemme $3.12$ de
\cite{9}.
\end{proof}

Soit $(G,\sigma)$ un fibr\'e quadratique, de spectre $r_{\ell}$.
Alors
$\Ext^{1}_{asym}(G,G)$ s'identifie \`a
$H^{1}(\Lambda^{2}G^{*}\otimes\omega_{_{\proj_{1}}})$. On obtient
$$\begin{diagram}
\ext^{1}_{asym}(G,G)&=&\displaystyle
h^{1}(\Lambda^{2}G^{*}\otimes\omega_{_{\proj_{1}}})\hfill\\
&=&\displaystyle
h^{1}\left(\osum_{\ell}\frac{r_{\ell}(r_{\ell}-1)}{2}
{\cal{O}}_{\proj_{1}}(-2\ell-2)
\osum
(\osum_{\ell<m}r_{\ell}r_{m}{\cal{O}}_{\proj_{1}}(-\ell-m-2))\right)\hfill\\
&=&\displaystyle
\sum_{\ell\ge0}\frac{r_{\ell}(r_{\ell}-1)}{2}(2\ell+1)
+\sum\begin{Sb} \ell<m\\ \ell+m\ge0\end{Sb} r_{\ell}r_{m}(\ell+m+1)\hfill\\
\end{diagram}
$$

On d\'eduit de ce calcul, en raison de la connexit\'e du spectre d'un
fibr\'e de Higgs quadratique semi-stable, le lemme suivant:
\begin{lemme}\label{rigide}
Le fibr\'e vectoriel $G$ sous-jacent \`a un fibr\'e de Higgs quadratique de
rang
$d$ \`a la fois semi-stable et quasi-rigide  $(G,\varphi,\tau)$
s'identifie ou bien \`a $R_{0}=d{\cal{O}}_{\proj_{1}}(-1)$  ou bien \`a
$R_{1}=
{\cal{O}}_{\proj_{1}}\osum (d-2){\cal{O}}_{\proj_{1}}(-1) \osum
{\cal{O}}_{\proj_{1}}(-2)$ (si $r\ge 3$).
\end{lemme}

\subsection{Description de l'ouvert $U_{a}(d)$.}

Soit $a\in\planp$ et $\proj_{1}$ une droite du plan projectif ne passant
pas par ce point. Consid\'erons la projection de centre a sur la droite:
$$\pi:\planp-\{a\}\lra\proj_{1}$$

L'espace total du fibr\'e normal $N={\cal{O}}_{\droitep}(1)$
s'identifie \`a l'ouvert $\planp-\{a\}$.

On va d\'ecrire l'ouvert $U_{a}(d)$ des th\^eta-caract\'eristiques semi-stables
de degr\'e $d$ dont le support ne passe pas par $a$ en termes de fibr\'es de
Higgs quadratiques sur $\droitep$.

\begin{prop}\label{cat-equiv}
La cat\'egorie des
th\^eta-caract\'eristiques dont le support sch\'ematique ne passe pas par le
point $a$ est \'equivalente,
par image directe, \`a la cat\'egorie des $N$-fibr\'es
de Higgs quadratiques sur $\droitep$.
\end{prop}
\begin{proof} Rappelons d'abord l'\'equivalence des cat\'egories entre
les faisceaux purs de dimension $1$  dont le support ne passe pas
par le point
$a$ et les fibr\'es de Higgs sur $\droitep$
(\cf la proposition 3.10 de \cite{9}). La correspondance
se fait par image directe: si ${\cal{F}}$ est un faisceau pur de
dimension 1 dont le support $C=\supps({\cal{F}})$ ne passe pas par $a$ alors
l'image directe $G=\pi_{*}({\cal{F}})$ est localement libre (par puret\'e) de
rang $\mult({\cal{F}})$ et munie d'une structure de
$\pi_{*}({\cal{O}}_{N})$-module. La donn\'ee d'une telle structure \'equivaut
 maintenant \`a se donner un morphisme $\varphi:G\lra G\otimes N$ d'o\`u la
structure de Higgs sur $G$. R\'eciproquement la donn\'ee d'un
fibr\'e vectoriel sur $\droitep$ muni d'une structure de
$\pi_{*}({\cal{O}}_{N})$-module d\'efinit un faisceau coh\'erent pur
${\cal{F}}$ sur N dont le support est fini au-dessus de $\droitep$: Soit
$\lambda$ la section canonique de $\pi^{*}(N)$ et consid\'erons le morphisme
$$\varphi-\lambda id_{G}:\pi^{*}(G)\lra\pi^{*}(G)\otimes\pi^{*}(N).$$
Le support sch\'ematique de ${\cal{F}}$ est alors d\'efini par le d\'eterminant
de ce morphisme (c'est la courbe spectrale associ\'e \`a $\varphi$) et
${\cal{F}}$ s'identifie au conoyau du morphisme
$(\varphi-\lambda id_{G})\otimes id_{N^{-1}}$.

Maintenant, le $\omega_{_{\planp}}$-dual de ${\cal{F}}$ s'identifie au
$\omega_{_{\droitep}}$-dual du fibr\'e de Higgs associ\'e:

\begin{lemme} Soit ${\cal{F}}$ un faisceau pur dont
le support ne passe pas par $a$ et $(G,\varphi)$ le fibr\'e de Higgs associ\'e
par la correspondance de la proposition \ref{cat-equiv}. Alors le fibr\'e de
Higgs associ\'e \`a ${\cal{F}}^{\vee}$ s'identifie au
$\omega_{_{\droitep}}$-dual $(G^{\vee},^{\vee}\varphi)$ de $(G,\varphi)$
\end{lemme}
\begin{proof} Soit $C$ le support sch\'ematique de ${\cal{F}}$.
Consid\'er\'e sur $C$, le faisceau ${\cal{F}}^{\vee}$ s'identifie \`a
$\ul{\Hom}_{{\cal{O}}_{C}}({\cal{F}},\omega_{C})$. On va consid\'erer le
morphisme fini $C\lra\droitep$ d\'efini par $\pi$. La dualit\'e de
Serre-Grothendieck fournit un isomorphisme
$$R\pi_{*}R\ul{\Hom}_{{\cal{O}}_{C}}({\cal{F}},\omega_{C})\simeq
R\ul{\Hom}_{{\cal{O}}_{\droitep}}(R\pi_{*}{\cal{F}},\omega_{\droitep})$$ dans
la cat\'egorie d\'eriv\'ee. Maintenant, comme ${\cal{F}}$ est un
${\cal{O}}_{C}$-module de Cohen-Macaulay, on a
$\ul{\Ext}^{i}({\cal{F}},\omega_{C})=0$ pour $i\ge 1$ \cite{12}. Par
cons\'equent,
les images
directes sup\'erieures \'etant nulles et $\pi_{*}{\cal{F}}$ \'etant localement
libre, l'isomorphisme ci-dessus se lit
$\pi_{*}\ul{\Hom}_{{\cal{O}}_{C}}({\cal{F}},\omega_{C})\simeq
\ul{\Hom}_{{\cal{O}}_{\droitep}}(\pi_{*}{\cal{F}},\omega_{\droitep})$.
Et puisque un isomorphisme dans la cat\'egorie d\'eriv\'ee entre deux
complexes de longueur $1$ d\'efinit un isomorphisme entre les objets
on a l'isomorphisme cherch\'e.
\end{proof}
La correspondance de la proposition est donc obtenu en associant \`a la
th\^eta-caract\'eristique $({\cal{F}},\sigma)$ le fibr\'e de Higgs quadratique
$(G,\varphi,\tau)$ o\`u $G=\pi_{*}({\cal{F}})$ et $\tau=\pi_{*}(\sigma)$.
\end{proof}

Remarquons que cette correspondance pr\'eserve la multiplicit\'e et la
caract\'eristique d'Euler-Poincar\'e (suite spectrale de Leray). En
particulier
${\cal{F}}$ est semi-stable (resp. stable) si et seulement si le
fibr\'e de Higgs associ\'e
$(G,\varphi)$ est semi-stable (resp. stable). De plus, si
$({\cal{F}},\sigma)$ est une th\^eta-caract\'eristique et $(G,\varphi,\tau)$
le fibr\'e de Higgs quadratique associ\'e alors si ${\cal{F}}'\subset{\cal{F}}$
est un sous-faisceau, le sous-faisceau de Higgs $\pi_{*}({\cal{F}}')$ de
$(G,\varphi)$ est d'orthogonal $\pi_{*}({\cal{F}}^{\prime\ort})$. En
particulier
$({\cal{F}},\sigma)$ est semi-stable (resp. stable) si et seulement si
$(G,\varphi,\tau)$ est semi-stable (resp. stable).

\subsubsection{Le foncteur $\protect\underline{HQ}(d,N)$}

Soient $d$ un entier et soit $N$ choisi de fa\c con \`a
ce que pour tout fibr\'e
de Higgs quadratique semi-stable $(G,\varphi,\tau)$ de rang $d$, $G(N)$
n'ait pas de cohomologie sup\'erieure et soit engendr\'e par ses sections
globales.
Ceci est possible, la familles des fibr\'es de Higgs quadratiques
semi-stables de rang fix\'e \'etant limit\'ee en raison de la
correspondance ci-dessus et le fait que la famille des
th\^eta-caract\'eristiques de multiplicit\'e fix\'ee est limit\'ee. Soit
$n=dN$,
$H=k^{n}$ et ${\tr{H}}=H\otimes_{k}{\cal{O}}_{\droitep}(-N)$. Consid\'erons,
pour toute vari\'et\'e alg\'ebrique $S$, les quadruplets
$(G,\alpha,\varphi,\tau)$ form\'es d'un quotient coh\'erent
${\tr{H}}_{S}\hfl{\alpha}{}G$, $S$-plat, d'une structure de Higgs
$\varphi$ sur $G$ et d'un isomorphisme de Higgs sym\'etrique
$\tau:G\lra G^{\vee}$ satisfaisant aux conditions suivantes:
pour tout point ferm\'e $s$ de $S$, $(G_{s},\varphi_{s},\tau_{s})$ est un
fibr\'e de Higgs quadratique de rang $d$ et $\alpha$ induit
un isomorphisme $H\otimes{\cal{O}}_{S}\simeq pr_{1*}(G(N))$. Deux
quadruplets  $(G,\alpha,\varphi,\tau)$ et
$(G',\alpha',\varphi',\tau')$ sont dits \'equivalents s'il existe un
isomorphisme de fibr\'es de Higgs
$f:(G,\varphi)\lra (G',\varphi')$ tel que
$\alpha'=f\circ \alpha$ et tel que
$\tau'\circ f=(^{\vee}f)^{-1}\circ\tau$. Notons
par $[G,\alpha,\varphi,\tau]$  la classe d'\'equivalence du
quadruplet $(G,\alpha,\varphi,\tau)$ et par
$\ul{HQ}(d,N)(S)$ l'ensemble des classes d'\'equivalence de
tels quadruplets. Ceci d\'efinit un foncteur:
$$ \ul{HQ}(d,N):\text{Vari\'et\'es alg\'ebriques}
\lra\text{Ensembles}.
$$
Le foncteur qu'on obtient en supposant de plus que pour tout point ferm\'e
$s$ on a $[G_{s},\varphi_{s},\tau_{s}]$ semi-stable sera not\'e
$\ul{HQ}^{ss}(d,N)$, celui qu'on obtient en consid\'erant
au lieu des quadruplets seulement les  triplets $[G,\alpha,\tau]$ (sans
structure de Higgs) sera not\'e $\ul{Q}(d,N)$.

On va repr\'esenter $\ul{HQ}^{ss}(d,N)$ de la mani\`ere suivante:

Soit $\Groth(d,N)$ le sch\'ema de
Hilbert-Grothendieck des quotients de $H\otimes{\cal{O}}_{\droitep}(-N)$ de
polyn\^ome de Hilbert $md$. Notons $\Groth_{0}$ l'ouvert
correspondant aux quotients $G$ tels que $H\simeq\H^{0}(G(N))$.
Cet ouvert param\`etre les classes d'\'equivalence des
 fibr\'es vectoriels $G$ de rang $d$ et de caract\'eristique
d'Euler-Poincar\'e nulle muni d'un isomorphisme $\alpha:H\simeq G(N)$.
Ici, $(G,\alpha)$ et $(G',\alpha')$ sont \'equivalents s'il existe un
isomorphisme $f:G\lra G'$ tel que $\alpha'=f\circ\alpha$.

Soit $\underline{G}$ le quotient universel sur
$\Groth_{0}\times\droitep$ et consid\'erons le faisceau coh\'erent
$pr_{1*}(\underline{G}^{\vee}(N))$. Par hypoth\`ese ce faisceau est
localement libre de fibre $H^{0}(G^{\vee}(N))$ au
dessus du point repr\'esent\'e par $[G,\alpha]$.
Soit ${\cal{R}}$ son fibr\'e de
rep\`eres.
L'espace total de ${\cal{R}}$ param\`etre les classes d'\'equivalence
des fibr\'es vectoriels $G$ de rang $d$ et de caract\'eristique
d'Euler-Poincar\'e nulle muni d'un isomorphisme $\alpha:H\simeq G(N)$
 et d'un isomorphisme $\beta:H\simeq G^{\vee}(N)$. Deux triplets sont
\'equivalents si s'il existe un isomorphisme $f:G\lra G'$ tel que
$\alpha'=f\circ\alpha$ et $\beta'=^{\vee}f^{-1}\circ\beta$.

On a une involution sur ${\cal{R}}$ en associant au triplet
$[G,\alpha,\beta]$ le triplet $[G^{\vee},\beta,\alpha]$. Soit
$Q$ le sch\'ema
des points fixes de ${\cal{R}}$ sous cette involution.

\begin{lemme}
Le sch\'ema $Q$ repr\'esente le foncteur $\ul{Q}$.
\end{lemme}
\begin{proof} Analogue \`a la d\'emonstration de la proposition 7.2 de
\cite{11}.
\end{proof}
Consid\'erons le triplet
universel $[\ul{G},\ul{\alpha},\ul{\tau}]$ sur
$$Q\times \droitep.$$
On peut identifier via $\ul\tau$ les faisceaux des
${\cal{O}}_{Q\times \droitep}$-modules
$\ul{\Hom}(\ul{G},\ul{G}(-3))$ et  $\ul{\Hom}(\ul{G},\ul{G}^{\vee}(-3))$.
Soit $\ul{\Hom}_{sym}(\ul{G},\ul{G}^{\vee}(-3))$ le
sous-faisceau des homomorphismes sym\'etriques. Ce sous-faisceau d\'efinit via
l'identification ci-dessus un sous-faisceau de
$\ul{\Hom}(\ul{G},\ul{G}(-3))$ qu'on note
$\ul{\Hom}_{sym}(\ul{G},\ul{G}(-3))$.
Soit $HQ(d,N)$ le fibr\'e
vectoriel (au sens de Grothendieck) associ\'e au faisceau coh\'erent
$R^{1}pr_{1*}\ul{\Hom}_{sym}(\ul{G},\ul{G}(-3))$.
Sa fibre au-dessus de $[G,\alpha,\tau]$ s'identifie, par dualit\'e
de Serre, aux morphismes $\varphi:G\lra G(1)$ $\tau$-sym\'etriques,
\ie aux morphismes $\varphi$ tels que $(G,\varphi,\tau)$ soit un
fibr\'e de Higgs quadratique.
Le sch\'ema $HQ(d,N)$ repr\'esente
$\ul{HQ}(d,N)$. On note $HQ^{ss}(d,N)$ l'ouvert de $HQ(d,N)$
repr\'esentant
$\ul{HQ}^{ss}(d,N)$.
Soit $\ul{T}_{a}^{ss}(d,N)$ le sous-foncteur ouvert de
$\ul{T}^{ss}(d,N)$ d\'efini par les triplets
$[{\cal{F}},\alpha,\sigma]$ tels que
$a\not\in\supps({\cal{F}})$.
D'apr\`es le proposition \ref{cat-equiv}, on a

\begin{prop}\label{isom} L'image directe induit un
isomorphisme de foncteurs entre les foncteurs $\ul{T}_{a}^{ss}(d,N)$ et
$\ul{HQ}^{ss}(d,N)$.
\end{prop}

Le groupe $GL(H)$ op\`ere sur $HQ^{ss}(d,N)$. D'apr\`es la
proposition \ref{isom}, le quotient de Mumford s'identifie \`a
l'ouvert
$U_{a}(d)$. On note $HQ^{os}(d,N)$ (resp. $U_{a}^{os}(d)$)
l'ouvert correspondant aux th\^eta-caract\'eristiques
${\cal{O}}$-stables. L'op\'eration de $GL(H)/\{\pm 1\}$ est libre sur
$HQ^{os}(d,N)$ de quotient
$U_{a}^{os}(d)$. Ainsi, le morphisme
$$HQ^{os}(d,N)\lra U_{a}^{os}(d)$$ est lisse.
On note $Q^{os}(d,N)$ l'ouvert de $Q(d,N)$ d\'efini par l'image de
$HQ^{os}(d,N)$ sous la projection.

\subsection{D\'emonstration du th\'eor\`eme \protect\ref{codimU}}

L'ouvert $U_{a}(d)$ \'etant \'equidimensionnel (car $\Theta(d)$ l'est)
et le ferm\'e des th\^eta-caract\'eristiques semi-stables, non
${\cal{O}}$-stables \'etant de codimension au moins $2$, il suffit de
consid\'erer l'ouvert
$U_{a}^{os}(d)$. Dans ce qui suit on notera simplement
$HQ^{os}$ et $Q^{os}$ les sch\'emas $HQ^{os}(d,N)$ et
$Q^{os}(d,N)$.

Consid\'erons la projection
$$p:HQ^{os}\lra Q^{os}$$
et un triplet $[G,\alpha,\tau]$ correspondant \`a un point ferm\'e de
$Q^{os}$. Soit $O(G,\alpha,\tau)$ l'orbite dans $Q^{os}$ de
$[G,\alpha,\tau]$ sous l'action de $GL(H)$ . Remarquons que $\tau$ est
n\'ecessairement donn\'e par une constante non nulle \ie $\tau\in k^{*}$.
En effet, $Q^{os}$ est par d\'efinition l'image de $HQ^{os}$, qui repr\'esente
les th\^eta-caract\'eristique ${\cal{O}}$-stables ne passant par par $a$. Or
si $({\cal{F}},\sigma)$ est ${\cal{O}}$-stable $\sigma\in k^{*}$.
Il s'ensuit que l'orbite de d\'epend que de $G$.
La codimension de l'orbite de $[G,\alpha,\tau]$ est
donn\'ee par $\ext^{1}_{asym}(G,G)$.
Soit $Z(G)$ l'image
r\'eciproque de $O(G,\alpha,\tau)$ sous $p$ et
calculons sa codimension: au-dessus de l'orbite la dimension de la fibre
est donn\'ee par $\hom_{sym}(G,G(1))$. D'apr\`es la proposition
\ref{isom} et \cite{11}, la vari\'et\'e
$HQ^{os}$ est lisse de  dimension  $\frac{d(d+3)}{2}+n^{2}$. La
vari\'et\'e
$Q^{os}$ est lisse, elle aussi: c'est un ouvert d'un sch\'ema des points
fixes sous une involution d'une vari\'et\'e lisse. Sa dimension se calcule
suivant la suite exacte suivante: (\cite{11} corollaire 7.6)
$$0\ra\Hom_{sym}(G,G)\ra
T_{[G,\alpha,\tau]}Q\ra
T_{[G,\alpha]}\Groth\ra
\Ext^{1}_{sym}(G,G)\ra 0$$
On en d\'eduit que $\dim Q^{os}=n^{2}-d^{2}+\chi_{sym}(G,G)$.
Par la formule de Riemann-Roch $\chi_{sym}(G,G)=\frac{d(d+1)}{2}$ d'o\`u
$\dim Q^{os}=n^{2}-\frac{d(d-1)}{2}$. De ces calculs on obtient
que la codimension de l'image r\'eciproque est donn\'ee par
$$\text{ext}^{1}_{asym}(G,G)-\text{ext}^{1}_{sym}(G,G(1)).$$
Soit $r_{\ell}$ le spectre de $G$. Ce spectre est connexe par
semi-stabilit\'e. L'espace vectoriel $\Ext^{1}_{asym}(G,G)$
s'identifie \`a $H^{1}(\Lambda^{2}G^{*}\otimes\omega_{_{\proj_{1}}})$ et
l'espace vectoriel $\Ext^{1}_{sym}(G,G(1))$ \`a
$H^{1}(S^{2}G^{*}\otimes\omega_{_{\proj_{1}}}(1))$. Par cons\'equent:
$$
\text{ext}^{1}_{asym}(G,G)-\text{ext}^{1}_{sym}(G,G(1))=
\displaystyle
\sum_{\ell\ge 0}\left(\frac{r_{\ell}(r_{\ell}-1)}{2}-2\ell r_{\ell}\right)+
\sum\begin{Sb}\ell+m\ge 0\\ \ell<m\end{Sb}r_{\ell}r_{m}.
$$
Cette expression est toujours positive,
nulle dans exactement $3$ cas:
celui o\`u $G$ est \'egale \`a $R_{0}$ ou $R_{1}$
du lemme \ref{rigide} ou si
$$G={\cal{O}}(\ell)\osum{\cal{O}}(\ell-1)\osum\dots
\osum{\cal{O}}(-\ell-1)\osum{\cal{O}}(-\ell-2).$$
Les images r\'eciproques $Z(G)$ sont, dans ces trois cas, $GL(H)/\{\pm
id\}$-invariants et irr\'eductibles. Et comme ils d\'efinissent respectivement
les th\^eta-caract\'eristiques ineffectives, semi-ineffectives et canoniques on
d\'eduit la premi\`ere partie du th\'eor\`eme, la codimension du
compl\'ementaire de
$U_{a}^{2,os}(d)$ \'etant, sauf dans le cas des
th\^eta-caract\'eristiques canoniques,
sup\'erieure ou \'egale \`a $1$ d'apr\`es ce qui
pr\'ec\`ede. Maintenant l'expression ci-dessus vaut $1$ exactement,
toujours par connexit\'e du spectre, pour
\begin{list}{-}{}
\item $G=R_{2}=
          2{\cal{O}}\osum(d-4){\cal{O}}(-1)\osum 2{\cal{O}}(-2)
          \text{ ($d\ge 5$ pour avoir la connexit\'e})$\par et aussi pour
\item $G=R_{3}=
        {\cal{O}}(1)\osum{\cal{O}}\osum2{\cal{O}}(-1)\osum{\cal{O}}(-2)
        \osum{\cal{O}}(-3)
        \text{ (ce qui impose $d=6$})$
\end{list}
Les images r\'eciproques $Z(G)$ sont, dans ces deux
cas, $GL(H)/\{\pm id\}$-invariants et irr\'eductibles.
Ils d\'efinissent, respectivement, les sous-sch\'emas localement ferm\'es
$U_{a}^{1,os}(d)\backslash U_{a}^{2,os}(d)$ et
$U_{a}^{2,os}(6)\backslash U_{a}^{3,os}(6)$ d'o\`u le th\'eor\`eme.
\cqfd

La d\'emonstration montre que le nombre maximal de sections
possibles d'une th\^eta-caract\'eristique semi-stable de degr\'e $d$ est
donn\'ee, si $d$ est impair par la dimension de l'espace des sections
de (avec $\ell=\frac{d-3}{2}$)
$$G={\cal{O}}(\ell)\osum{\cal{O}}(\ell-1)\osum\dots
\osum{\cal{O}}\osum{\cal{O}}(-1)\osum{\cal{O}}(-2)\dots
\osum{\cal{O}}(-\ell-1)\osum{\cal{O}}(-\ell-2),$$
si $d$ est pair par la dimension de l'espace des sections de
$$G={\cal{O}}(\ell)\osum{\cal{O}}(\ell-1)\osum\dots
\osum{\cal{O}}\osum2{\cal{O}}(-1)\osum{\cal{O}}(-2)\dots
\osum{\cal{O}}(-\ell-1)\osum{\cal{O}}(-\ell-2),$$
On en d\'eduit le corollaire suivant:
\begin{cor}
Soit $({\cal{F}},\sigma)$ une th\^eta-caract\'eristique semi-stable de degr\'e
$d$.
\begin{list}{-}{}
\item Si $d$ est impair on a
$\h^{0}({\cal{F}})\le\frac{(d-3)^2}{8}+\frac{3(d-3)}{4}+1$.
\item Si $d$ est pair on a
$\h^{0}({\cal{F}})\le\frac{(d-4)^2}{8}+\frac{3(d-4)}{4}+1$.
\end{list}
\end{cor}

En particulier, $\Theta^{1}(d)$ est vide si $d\leq 4$,
$\Theta^{2}(d)$ est vide si
$d\leq 5$ et $\Theta^{3}(d)$ est vide si $d\leq 7$.

\begin{lemme}
Si $d\ge 5$ il existe une th\^eta-caract\'eristique
$({\cal{F}},\sigma)$ semi-stable de degr\'e
$d$ telle que $\h^{0}({\cal{F}})=2$.
\end{lemme}
\begin{proof}
Si $d\ge 6$, on peut consid\'erer par exemple la somme directe orthogonale
${\cal{O}}_{C}\osum{\cal{O}}_{C}\osum(d-6){\cal{O}}_{\ell}(-1)$ avec
$C$ une courbe de degr\'e $3$ et $\ell$ une droite de $\planp$. Si $d=5$,
consid\'erons le fibr\'e vectoriel $G=2{\cal{O}}\osum{\cal{O}}(-1)\osum
2{\cal{O}}(-2)$ sur $\droitep$ et l'isomorphisme
$\sigma:G\lra G^{\vee}$ d\'efini par $1$ sur l'anti-diagonale, $0$
sinon. Soit de plus $\varphi:G\lra G(1)$ d\'efini par la matrice
suivante:
$$\left(\begin{matrix}
\alpha&0&0&0&0\\
0&\beta&0&0&0\\
\lambda&\mu&0&0&0\\
0&0&\mu&0&0\\
0&0&\lambda&0&0
\end{matrix}\right)
$$
avec $\lambda$ et $\mu$ des constantes non nuls, $\alpha$ et
$\beta$ non proportionnels. On v\'erifie cas par cas qu'aucun
sous-fibr\'e de $G$ contenant ${\cal{O}}$ n'est de Higgs. Ainsi
$(G,\varphi,\sigma)$ est un fibr\'e de Higgs quadratique semi-stable,
d\'efinissant, pour
$a\in\planp$ fix\'e, par le proc\'ed\'e ci-dessus une
th\^eta-caract\'eristique semi-stable de degr\'e
$5$ ne passant pas par $a$.
\end{proof}

\subsection{Le compl\'ementaire de $U_{a}(d)$ dans $\Theta(d)$.}

Consid\'erons le compl\'ementaire $F_{a,p}(d)$  de  $U_{a}(d)$
dans $\Theta_{p}(d)$ (resp.
$F_{a,i}(d)$ de $U^{a}_{1}(d)$ dans $\Theta_{i}(d)$).

\begin{prop} Pour les ferm\'es $F_{a,p}(d)$ et $F_{a,i}(d)$ on a:
\begin{list}{\arabic{lc}.}{\usecounter{lc}}
\item le ferm\'e $F_{a,p}(d)$ est une hypersurface irr\'eductible de
$\Theta_{p}(d)$.
\item le ferm\'e $F_{a,i}(d)$ est une hypersurface irr\'eductible de
$\Theta_{i}(d)$.
\end{list}
\end{prop}
\begin{proof}
$1)$ Soient $a,b\in\planp$. Quitte \`a \'echanger $a$ et $b$ il s'agit de
montrer l'\'enonc\'e ci-dessus pour les th\^eta-caract\'eristiques
ne passant pas par $b$. D'abord on a besoin du lemme suivant:

 \begin{lemme} On a, pour $[R_{0},\alpha,\tau]\in Q^{os}$,
l'isomorphisme suivant, avec $K=\Aut(R_{0},\tau)$:
$$Sym^{os}(R_{0},R_{0}(1))/K\simeq I_{a}^{os}(d)$$
o\`u $I_{a}^{os}(d)$ d\'esigne l'ouvert de $U_{a,p}(d)$ correspondant
aux th\^eta-caract\'eristiques ineffectives.
\end{lemme}
\begin{proof}
 Le groupe $K$ s'identifie au stabilisateur de $[R_{0},\alpha,\tau]$
sous l'action de $GL(H)$ sur $Q^{os}$. Ce groupe agit sur
$Sym^{os}(R_{0},R_{0}(1))$ par conjugaison et l'on a le diagramme commutatif
$$
\begin{diagram}
Sym^{os}(R_{0},R_{0}(1))&\efl{i}{}&Z(G)\\
&\sefl{\pi'}{}&\sfl{\pi}{}\\
&&I_{a}^{os}(d)\\
\end{diagram}
$$
L'inclusion $i$ est compatible aux actions de $K$ et $GL(H)$ et $\pi$ est un
quotient g\'eom\'etrique.
Par suite, $\pi'$ est un quotient g\'eom\'etrique aussi,
d'o\`u le lemme.
\end{proof}
Si $(R,\varphi,\tau)$ est le fibr\'e de Higgs quadratique associ\'e \`a la
th\^eta-caract\'eristique $({\cal{F}},\sigma)$, le support du faisceau
${\cal{F}}$ est d\'efini par les couples $(x,\lambda)$ tels que
$$\det(\varphi(x)-\lambda\id_{R(x)})=0.$$ Maintenant
le fibr\'e
$\ul{Sym}(R_{0},R_{0}(1))$ est engendr\'e par ses sections globales.
Par cons\'equent l'image r\'eciproque sous le morphisme \'evaluation
$$Sym(R_{0},R_{0}(1))\lra \ul{Sym}(R_{0},R_{0}(1))(b)$$
de l'hypersurface  des applications lin\'eaires sym\'etriques de $R_{0}\lra
R_{0}(1)$ de d\'eterminant nul est une hypersurface irr\'eductible de
$Sym^{os}(R_{0},R_{0}(1))$. Cette hypersurface est $K$-invariante et
d\'efinit, par passage au quotient, l'hypersurface des
th\^eta-caract\'eristiques de
$I_{a}^{os}(d)$ dont le support passe par $b$. Cette
hypersurface est donc irr\'eductible et comme le ferm\'e des
th\^eta-caract\'eristiques non
${\cal{O}}$-stables est de codimension au moins deux
ceci est encore vrai sans supposer la ${\cal{O}}$-stabilit\'e.
Pour voir que l'hypersurface
$F_{b}(d)$  est encore irr\'eductible dans $\Theta_{p}$, il suffit de
montrer qu'il existe des th\^eta-caract\'eristiques ayant au moins deux
sections lin\'eairement ind\'ependantes ne passant ni par $a$ ni par $b$. Pour
cela, consid\'erons $Sym(R_{2},R_{2}(1))$. Ce fibr\'e n'est pas
engendr\'e par ses sections. Cependant l'image de
$$Sym(R_{2},R_{2}(1))\lra \ul{Sym}(R_{2},R_{2}(1))(b)$$
est un espace lin\'eaire de codimension $3$ qui n'est pas inclus dans
l'hypersurface  des applications lin\'eaires sym\'etriques de $R_{2}\lra
R_{2}(1)$ de d\'eterminant nul. Ceci d\'emontre la derni\`ere assertion
et par cons\'equent $(i)$.

$2)$ Un argument analogue \`a celui ci-dessus montre qu'on a, pour
$[R_{1},\alpha,\tau]\in Q^{os}$, l'isomorphisme suivant, avec
$K=\Aut(R_{1},\tau)$: $$Sym^{os}(R_{1},R_{1}(1))/K\simeq
SI_{a}^{os}(d),$$
o\`u $SI_{a}^{os}(d)$ d\'esigne l'ouvert de $U_{a,i}(d)$ correspondant
aux th\^eta-caract\'eristiques semi-ineffectives.

Maintenant le fibre
$\ul{Sym}(R_{1},R_{1}(1))$ n'est plus engendr\'e par ses sections.
Cependant l'intersection de l'hyperplan image de
$Sym(R_{1},R_{1}(1))$ dans $\ul{Sym}(R_{1},R_{1}(1))(b)$
avec
l'hypersurface
des application lin\'eaires sym\'etriques de $R_{1}\lra R_{1}(1)$ de
d\'eterminant nul est irr\'eductible de codimension $1$ dans l'image.
L'image r\'eciproque de cette intersection d\'efinit une hypersurface
de
$Sym^{os}(R_{1},R_{1}(1))$.
Cette hypersurface est $K$-invariante et
d\'efinit l'hypersurface des th\^eta-caract\'eristiques
de $SI_{a}^{os}(d)$ passant par $b$. Cette hypersurface est donc
irr\'eductible et ceci est encore vrai sans supposer la
${\cal{O}}$-stabilit\'e. Si $d\not=6$ ceci d\'emontre
$(ii)$, si
$d=6$ il faut encore voir  qu'il existe des
th\^eta-caract\'eristiques ayant $3$ sections lin. ind\'ependantes ne
passant ni par $a$ ni par $b$.
Pour cela, consid\'erons l'\'evaluation
$$Sym(R_{3},R_{3}(1))\lra \ul{Sym}(R_{3},R_{3}(1))(b)$$
L'image est un espace lin\'eaire de codimension $5$ n'\'etant pas
inclus dans l'hypersurface  des application lin\'eaires sym\'etriques
de $R_{2}\lra R_{2}(1)$ de d\'eterminant nul. Par cons\'equent il
existe de telles th\^eta-caract\'eristiques dans $U_{a}(6)$ qui
sont dans
$U_{b}(6)$, d'o\`u $(ii)$.
\end{proof}
\begin{cor} On a
\begin{list}{\arabic{lc})}{\usecounter{lc}}
\item L'hypersurface $\Theta^{1}(d)$ de $\Theta_{p}(d)$
est irr\'eductible, le ferm\'e $\Theta^{3}(d)$ est de codimension au
moins deux.
\item Le ferm\'e $\Theta^{2}(d)$ est de codimension au moins deux, sauf si
$d=6$ o\`u la codimension est $1$.
\end{list}
\end{cor}
\begin{proof}
L'\'enonc\'e est vrai pour l'ouvert $U_{a,p}(d)$. Il suffit donc de
montrer que
$\Theta^{1}(d)$ coupe $F_{a,p}(d)$ suivant un ferm\'e de codimension au
moins $2$. Soit $\ell$ une droite de $\planp$ passant par $a$,
$({\cal{F}},\sigma)$ la th\^eta-caract\'eristique somme directe orthogonale
$d{\cal{O}}(-1)$. Alors  $({\cal{F}},\sigma)\in F_{a,p}(d)$. Par
semi-continuit\'e $F_{a}^{1}(d)$ est un ferm\'e de $F_{a,p}(d)$, par
irr\'eductibilit\'e il est de codimension au moins $1$ dans $F_{a,p}(d)$.
L'assertion sur $\Theta^{3}(d)$ se d\'eduit de l'assertion analogue pour
l'ouvert $U_{a,p}(d)$ et de ce qui pr\'ec\`ede.
La d\'emonstration de $2)$ est analogue \`a celle de $1)$.
\end{proof}

\np
\section{La
r\'esolution minimale d'une th\^eta-caract\'eristique semi-stable}
\label{resolution}

Soit $({\cal{F}},\sigma)$ une th\^eta-caract\'eristique de degr\'e $d$.
Consid\'erons la $k$-alg\`ebre gradu\'e $S=k[X_{0},X_{1},X_{2}]$ et le
$S$-module
gradu\'e $M=\osum_{n\in\reln}\H^{0}({\cal{F}}(n))$. La r\'esolution minimale
gradu\'e de $M$ fournit une r\'esolution $(\star)$:
$$
0\lra E \hfl{\phi}{}
F\lra {\cal{F}}\lra 0
$$
avec $E=\osum_{i=1}^{n}{\cal{O}}_{\planp}(-a_{i})$ et
$F=\osum_{i=1}^{n+1}{\cal{O}}_{\planp}(-b_{i})$. Le morphisme $\phi$
s'identifie \`a une $n\times n$-matrice $\phi_{ij}$ avec
$\phi_{ij}\in\Hom({\cal{O}}(-a_{i}),{\cal{O}}(-b_{i}))$. Par
minimalit\'e de la r\'esolution, il ne peut exister de $i,j$ tel que
$\phi_{ij}$ soit une constante non nulle. On utilisera ce fait un peu
plus loin.

En appliquant le
foncteur $\ul{\Hom}(\cdot,\omega_{_{\planp}})$ \`a la
r\'esolution $(\star)$ on voit,
comme deux r\'esolution minimales sont isomorphes, que $F\simeq E^{\vee}$, \ie
que
$-b_{i}=a_{i}-3$ apr\`es r\'enumerotation. Posons
$b=h^{1}({\cal{F}})=h^{0}({\cal{F}})$,
$c=h^{1}({\cal{F}}(1))=h^{0}({\cal{F}}(-1))$.
On a $b>c$.
\begin{prop} Supposons $({\cal{F}},\sigma)$ semi-stable.
Alors la r\'esolution minimale $(\star)$ est de la forme suivante:
\begin{list}{-}{}
\item si $b=0$ on a $n=d$ et $(a_{i})=(2,\dots,2)$
\item si $b=1$ on a $n=d-2$ et $(a_{i})=(2,\dots,2,3)$
\end{list}
De plus, on peut supposer $\phi$ sym\'etrique dans le deux cas.
\end{prop}
\begin{proof}
Cette proposition est due \`a Dixon, si $b=0$ \cite{3}, \`a
Catanese \cite{2} sinon (\cf aussi Laszlo
\cite{7}) dans le cas o\`u le support sch\'ematique de
$({\cal{F}},\sigma)$ est lisse. L'argument de Laszlo se g\'en\'eralise
\`a notre situation, quitte \`a remplacer
l'argument d'irr\'eductibilit\'e du support qu'il invoque par un
argument utilisant la semi-stabilit\'e de $({\cal{F}},\sigma)$.
Nous omettons les d\'etails ici en raison de la proposition suivante.
\end{proof}

\begin{prop}\label{Beilinson}
Soit $({\cal{F}},\sigma)$ une th\^eta-caract\'eristique semi-stable.
Si $({\cal{F}},\sigma)$ est ineffective, ${\cal{F}}$ est conoyau du
morphisme  $d_{2}$ sym\'etrique et g\'en\'eriquement injectif,
donn\'e par la suite spectrale de Beilinson et $\sigma$:
$$\H^{1}({\cal{F}}(-1))\otimes{\cal{O}}_{\planp}(-2)
\hfl{d_{2}}{}\H^{1}({\cal{F}}(-1))^{*}\otimes {\cal{O}}_{\planp}(-1)$$\par
Si $({\cal{F}},\sigma)$ est semi-ineffective, alors ${\cal{F}}$ est
conoyau du morphisme $d_{2}$ sym\'etrique et g\'en\'eriquement injectif,
donn\'e par la suite  spectrale de Beilinson et $\sigma$:
$$K\otimes{\cal{O}}_{\planp}(-2)\osum
H^{1}({\cal{F}})\otimes{\cal{O}}_{\planp}(-3)
\hfl{d_{2}}{}
K^{*}\otimes{\cal{O}}_{\planp}(-1)\osum
H^{1}({\cal{F}})^{*}\otimes{\cal{O}}_{\planp},
$$ o\`u
$K$ est le noyau de
l'application canonique
$\H^{1}({\cal{F}}(-1))\lra\H^{1}({\cal{F}})\otimes V$ avec
$V=H^{0}(Q)$, $Q$ \'etant le quotient tautologique sur $\planp$.
\end{prop}
\begin{proof}
 Soit ${\cal{E}}$ un ${\cal{O}}_{\planp}$-module
coh\'erent. La suite  spectrale de Beilinson est donn\'ee, en degr\'e $1$,
par:
$$E_{1}^{p,q}=\H^{q}(\planp,{\cal{E}}(p))\otimes\Lambda^{-p}Q^{*}.$$
Elle est d'aboutissement
${\cal{E}}$ en degr\'e $0$, nul sinon.

Soit $({\cal{F}},\sigma)$ une th\^eta-caract\'eristique de degr\'e $d$.
Supposons que
$h^{0}(\planp,{\cal{F}})\leq1$. Par d\'ecroissance de la suite des
$h^{1}({\cal{F}}(i))$, on a
$h^{1}(\planp,{\cal{F}}(1))=0$. Nous allons appliquer la suite
spectrale de Beilinson \`a ${\cal{F}}(1)$, puis tordre par
${\cal{O}}_{\planp}(-1)$.
$$
\begin{array}{|c|c|c|}
\hline
0 & 0 & 0\\
\hline
\H^{1}({\cal{F}}(-1))\otimes\Lambda^{2}Q^{*}(-1)&
\H^{1}({\cal{F}})\otimes Q^{*}(-1)&0\\
\hline
0&\H^{0}({\cal{F}})\otimes Q^{*}(-1)&
\H^{0}({\cal{F}}(1))\otimes {\cal O}(-1)\\
\hline
\end{array}
$$
Soient $H_{1}$ et $H_{2}$ les
${\cal{O}}_{\planp}$-modules de cohomologie du complexe
$$
 0\lra \H^{1}({\cal{F}}(-1))\otimes\Lambda^{2}Q^{*}(-1)\hfl{d_{1}^{-2,1}}{}
       \H^{1}({\cal{F}})\otimes Q^{*}(-1)\lra 0,
$$ et $K_{2}$ et $K_{3}$ les ${\cal{O}}_{\planp}$-modules de
cohomologie du complexe
$$
 0\lra \H^{0}({\cal{F}})\otimes Q^{*}(-1)\hfl{d_{1}^{-1,0}}{}
       \H^{0}({\cal{F}}(1))\otimes {\cal{O}}_{\planp}(-1) \lra 0.
$$  Le terme $E_{2}^{p,q}$ s'\'ecrit par cons\'equent:
$$
\begin{array}{|c|c|c|}
\hline
0  &   0    &    0\\
\hline
H_{1}  &  H_{2}  &    0\\
\hline
0    &  K_{2}  & K_{3}\\
\hline
\end{array}
$$
Soient $L_{1}$ et $L_{3}$ les
${\cal{O}}_{\planp}$-modules de cohomologie du complexe
$$0 \lra H_{1}\hfl{d_{2}}{} K_{3} \lra 0$$ Par suite, on a pour le
terme $E_{3}$:
$$
\begin{array}{|c|c|c|}
\hline
0    &    0    &    0\\
\hline
L_{1}&  H_{2}  &    0\\
\hline
0    &  K_{2}  &    L_{3}\\
\hline
\end{array}
$$
Du fait que $E_{3}=E_{\infty}$, on d\'eduit que
$K_{2}=L_{1}=0$ et qu'on a la suite exacte
$$0\lra L_{3}\lra{\cal{F}}\lra H_{2}\lra 0.$$

D\'eterminons $K_{3}:$ Si $b=0$, on a bien s\^ur
$K_{3}=\H^{0}({\cal{F}}(1))\otimes{\cal{O}}_{\planp}(-1)$. Si $b=1$,
consid\'erons le diagramme commutatif avec lignes et colonnes exactes, la
premi\`ere ligne \'etant la dualis\'e tordue par
$H^{0}({\cal{F}})\otimes{\cal{O}}_{\planp}(-1)$ de la suite exacte
d'Euler:
$$
\begin{diagram}
 &  & & & 0& & 0\\
 &  & & & \vfl{}{}& & \vfl{}{}\\
 &  & & & M_{1}&\ra& M_{2}\\
 &  & & & \vfl{}{}& & \vfl{}{}\\
 0&\ra &\H^{0}({\cal{F}})\otimes Q^{*}(-1)
 &\ra&\H^{0}({\cal{F}})\otimes
 V^{*}\otimes {\cal{O}}(-1)&\ra&
 \H^{0}({\cal{F}})\otimes {\cal{O}} &\ra&0\\
 & & \Vert &&\vfl{}{}& &\vfl{}{}\\
 0&\ra &\H^{0}({\cal{F}})\otimes
 Q^{*}(-1)&\ra&\H^{0}({\cal{F}}(1))\otimes{\cal{O}}(-1)&
 \ra&K_{3}&\ra& 0\\
 & & & & \vfl{}{}& & \vfl{}{}\\
 & & & & N_{1}&\ra& N_{2}\\
 & & & & \vfl{}{}& & \vfl{}{}\\
 & & & & 0& & 0\\
\end{diagram}
$$

Par le lemme du serpent, $M_{1}\simeq M_{2}$ et
$N_{1}\simeq N_{2}$.
Supposons d'abord que  $M_{1}$ (et donc aussi $N_{1}$) nul. A ce moment-l\`a,
$N_{1}$ (et aussi $N_{2}$) est isomorphe \`a
$H\otimes{\cal{O}}_{\planp}(-1)$, o\`u $H$ est le conoyau de l'application
canonique $\H^{0}({\cal{F}})\otimes V^{*}\lra \H^{0}({\cal{F}}(1)).$
Par cons\'equent,
$K_{3}$ s'identifie \`a une extension de
$\H^{0}({\cal{F}})\otimes{\cal{O}}_{\planp}$ et
$H\otimes{\cal{O}}_{\planp}(-1)$. L'espace de telles extensions \'etant
nul, on a, si $b$=1,
$$K_{3}=\H^{0}({\cal{F}})\otimes{\cal{O}}_{\planp}\osum
H\otimes{\cal{O}}_{\planp}(-1).$$
Montrons que $M_{1}$ non nul est impossible. Sinon, ce faisceau est
n\'ecessairement isomorphe \`a
${\cal{O}}_{\planp}(-1)$. Mais ceci signifie que $K_{3}$ s'identifie
\`a ${\cal{O}}_{\ell}\osum{\cal{O}}(-1)^{d-2}$, $\ell$
\'etant une droite. Puisque $H_{1}$ est localement libre, ce
${\cal{O}}_{\ell}$ est contenu dans $L_{3}$ et par suite dans
${\cal{F}}$. Mais ceci est en contradiction avec la semi-stabilit\'e de
${\cal{F}}$, la caract\'eristique d'Euler-Poincar\'e de
${\cal{O}}_{\ell}$ \'etant strictement positive.

\begin{lemme}\label{Bdual}
La suite spectrale de Beilinson appliqu\'ee \`a ${\cal{F}}^{\vee}(1)$ puis
tordue
par ${\cal{O}}_{\planp}(-1)$ est naturellement isomorphe
 \`a la suivante, dont le terme $E_{1}$ est donn\'e par
\smallskip
$$
\begin{diagram}
\H^{0}({\cal{F}}(1))^{*}
\otimes\Lambda^{2}Q^{*}(-1)&\efl{^{\vee}d_{1,{\cal{F}}}^{-1,0}}{}&
\H^{0}({\cal{F}})^{*}\otimes Q^{*}(-1)&&0\\
0&&\H^{1}({\cal{F}})^{*}
\otimes Q^{*}(-1)&\efl{^{\vee}d_{1,{\cal{F}}}^{-2,1}}{}&
\H^{1}({\cal{F}}(1))^{*}\otimes {\cal O}(-1)\\
\end{diagram}
$$
et le terme $E_{2}$ par
$$\Ker(^{\vee}d_{1,{\cal{F}}}^{-1,0})
\efl{^{\vee}d_{2}^{-2,1}}{}\Coker(^{\vee}d_{1,{\cal{F}}}^{-2,1}).$$
\end{lemme}

Admettons pour un instant ce lemme est terminons la d\'emonstration de la
proposition.
Montrons d'abord que $\sigma$ d\'efinit un isomorphisme
de $K_{3}^{\vee}$ avec $H_{1}$:
par functorialit\'e, $\sigma$ fournit un isomorphisme au niveau des suites
spectrales de Beilinson associ\'es \`a ${\cal{F}}$ et ${\cal{F}}^{\vee}$ ce qui
donne, en utilisant le lemme pr\'ec\'edent, le diagramme commutatif

$$
\begin{diagram}
\H^{1}({\cal{F}}(-1))\otimes \Lambda^{2}Q^{*}(-1)&\efl{d^{-2,1}_{1}}{}&
\H^{1}({\cal{F}})\otimes Q^{*}(-1)\\
\sfl{}{}&&\sfl{}{}\\
\H^{0}({\cal{F}}(-1))^{*}\otimes \Lambda^{2}Q^{*}(-1)&
\efl{^{\vee}d_{1,{\cal{F}}}^{-1,0}}{}&
\H^{0}({\cal{F}})^{*}\otimes Q^{*}(-1)\\
\end{diagram}
$$

Le morphisme induit au niveau des noyaux, donne l'isomorphisme
$\tau:H_{1}\lra K_{3}^{\vee}$ cherch\'e. On en d\'eduit, en consid\'erant la
suite
exacte qui d\'efinit $H_{1}$ et $H_{2}$ que $H_{2}$ est de degr\'e nul et de
caract\'eristique d'Euler-Poincar\'e nulle. Comme $H_{2}$ est de torsion en
tant que conoyau de ${\cal{F}}$, il est n\'ecessairement nul.

Par cons\'equent, ${\cal{F}}$ s'identifie au conoyau du morphisme
$d_{2}^{-2,1}:H_{1}\lra K_{3}$. En utilisant la functorialit\'e des suites
spectrales de Beilinson et la sym\'etrie de $\sigma$, on obtient le diagramme
commutatif suivant:

$$
\begin{diagram}
H_{1}&\efl{d_{2}^{-2,1}}{}&K_{3}\\
\sfl{\tau}{}&\sefl{d_{2}}{}&\sfl{^{\vee}\tau}{}\\
K_{3}^{\vee}&\efl{}{^{\vee}d_{2}^{-2,1}}&H_{1}^{\vee}\\
\end{diagram}
$$
\bigskip
On en d\'eduit que ${\cal{F}}$ s'identifie au conoyau du morphisme sym\'etrique
$d_{2}$
\end{proof}

{\em D\'emonstration du lemme \ref{Bdual}:} Pour d\'emontrer ce lemme, revenons
\`a
la d\'efinition de la suite spectrale de Beilinson: soit $D^{\cdot}$
le complexe
$$0\lra\Lambda^{2}Q^{*}\etimes{\cal{O}}_{\planp}(-2)\lra
Q^{*}\etimes{\cal{O}}_{\planp}(-1)\lra{\cal{O}}_{\planp}
\etimes{\cal{O}}_{\planp}\lra 0$$
donn\'e par la r\'esolution de la diagonale $\Delta\subset\planp\times\planp$
de Beilinson.
Alors on a $R_{p_{1*}}(D^{\cdot}\etimes{\cal{F}})={\cal{F}}$ dans
${\cal{D}}^{b}_{c}(\planp)$ et la suite spectrale de Beilinson est donn\'ee
par la suite spectrale des foncteurs hyperderiv\'es. Dans notre cas, la suite
spectrale consid\'er\'ee est celle obtenue en regardant le complexe
$$L^{\cdot}=D^{\cdot}\otimes_{{\cal{O}}_{\planp\times\planp}}
({\cal{O}}(-1)\etimes{\cal{O}}(1))\etimes{\cal{F}}.$$ Ceci s'applique de
m\^eme au faisceau
${\cal{F}}^{\vee}$. Consid\'erons maintenant le complexe
$M^{\cdot}=R\ul{\Hom}(L^{\cdot},p_{1}^{!}(\Lambda^{2}Q^{*}(-2)))$
dans $D^{b}_{c}(\planp\times\planp)$ o\`u
$p_{1}^{!}(G)=p_{1}^{*}(G)\otimes p_{2}^{*}(\omega_{_{\planp}})[2]$ pour un
${\cal{O}}_{\planp}$-module coh\'erent $G$. Par dualit\'e de
Grothendieck-Serre,
on a, dans
$D^{b}_{c}(\planp)$,
$$Rp_{1*}(M^{\cdot})=R\ul{\Hom}(Rp_{1*}(L^{\cdot}))={\cal{F}}^{\vee}$$
De plus, dans
$D^{b}_{c}(\planp\times\planp)$, on a
$$M^{\cdot}\simeq\ul{\Hom}(D^{\cdot},\Lambda^{2}Q^{*}\etimes
{\cal{O}}_{\planp}(-2))\otimes_{{\cal{O}}_{\planp\times\planp}}
({\cal{O}}(-1)\etimes{\cal{O}}(1))\etimes
R\ul{\Hom}({\cal{F}},p_{2}^{*}(\omega_{_{\planp}}))[2].$$ La suite spectrale
des foncteurs hyperd\'eriv\'es du dernier complexe est celle
cherch\'ee et l'isomorphisme de suites spectrales du lemme se d\'eduit de
l'isomorphisme naturel
$$D^{\cdot}
\simeq \ul{\Hom}(D^{\cdot},\Lambda^{2}Q^{*}\etimes
{\cal{O}}_{\planp}(-2))[2],$$
d\'eduit de l'isomorphisme $Q^{*}\etimes{\cal{O}}(-1)\simeq
Q\etimes{\cal{O}}(+1)\otimes_{{\cal{O}}_{\planp\times\planp}}
\Lambda^{2}Q^{*}\etimes{\cal{O}}(-2).$
\cqfd

Dans
la suite, on aura besoin d'une description de la r\'esolution
minimale d'une th\^eta-caract\'eristique $({\cal{F}},\sigma)$
semi-stable dans le cas o\`u
$b=2$. A ce moment l\`a,
$h^{1}({\cal{F}}(2))=h^{0}({\cal{F}}(-2))=0$ et, a priori, deux cas se
pr\'esentent suivant si $c=0$ ou $c=1$.

Consid\'erons d'abord le cas o\`u $c=1$.
Remarquons qu'on a $-a_{i}+3\ge-1$, \ie $a_{i}\le4$ ici,
$h^{1}({\cal{F}}(2))$ \'etant nul. De plus, on a forc\'ement $a_{i}\ge0$. En
effet, s'il existait un
$i$ tel que $a_{i}\le-1$ on aurait $a_{i}+a_{j}-3\le a_{j}-4\le0$
pour tout $j$ ce qui est en contradiction \`a la minimalit\'e de la
r\'esolution. Le fibr\'e $E$ est donc de la forme
$\displaystyle \osum_{i=0}^{4}n_{i}{\cal{O}}_{\planp}(-i)$. En
tordant la r\'esolution par ${\cal{O}}_{\planp}(-1)$, on
voit que $n_{4}=1$ car $c=1$. De plus, comme $\h^{0}({\cal{F}})=2$ ceci
impose $n_{0}\ge 1$ et $n_{0}-n_{3}=1$. Encore par minimalit\'e de la
r\'esolution $n_{0}>1$ est impossible, d'o\`u n\'ecessairement $n_{0}=1$ et
$n_{3}=0$. On
a donc la factorisation suivante:
$$ \begin{diagram} && 0        &
& 0        & &\\ && \vfl{}{} & & \vfl{}{} & &\\
0&\lra& {\cal{O}}_{\planp}&\hfl{}{}&{\cal{O}}_{\planp}(1)&\hfl{}{}
& {\cal{O}}_{\ell}(1)&\lra&0\\
&& \vfl{}{} & & \vfl{}{} & & \vfl{}{}\\
0&\lra& E &\hfl{}{}&E^{\vee}&\hfl{}{}& {\cal{F}}&\lra&0\\
\end{diagram}
$$
Le noyau de ${\cal{O}}_{\ell}(1)\lra{\cal{F}}$ est, par le lemme
du serpent, dans le conoyau de ${\cal{O}}_{\planp}\lra E$. Par cons\'equent,
puisque ce conoyau est localement libre, ${\cal{O}}_{\ell}(1)$
s'injecte dans ${\cal{F}}$ ce qui est en contradiction avec la
semi-stabilit\'e de $({\cal{F}},\sigma)$.
Par suite, $c$ est n\'ecessairement \'egale \`a $0$. Dans ce
cas l\`a, on a $1\le a_{i}\le3$,
$\displaystyle E=\osum_{i=1}^{3}n_{i}{\cal{O}}_{\planp}(-i)$ avec
$n_{3}=2$ et la factorisation suivante: $$
\begin{diagram}
&& 0        & & 0        & &\\
&& \vfl{}{} & & \vfl{}{} & &\\
0&\lra&
n_{1}{\cal{O}}_{\planp}(-1)&\hfl{}{}&2{\cal{O}}_{\planp}&\hfl{}{} &
K&\lra&0\\
&& \vfl{}{} & & \vfl{}{} & & \vfl{}{}\\
0&\lra& E &\hfl{}{}&E^{\vee}&\hfl{}{}& {\cal{F}}&\lra&0\\
\end{diagram}
$$
Si $n_{1}=2$, alors on a $K\subset{\cal{F}}$ mais $\chi(K)>0$ ce qui
contredit la semi-stabilit\'e de ${\cal{F}}$. Il reste finalement le cas o\`u
$n_{1}=1$ ou $n_{1}=0$. On a donc d\'emontr\'e la proposition suivante:
\begin{prop}
Soit $({\cal{F}},\sigma)$ une th\^eta-caract\'eristique semi-stable telle que
$\h^{0}({\cal{F}})=2$. Alors la r\'esolution minimale de
${\cal{F}}$ qui est soit de type $(I)$
$$
0\ra (\osum_{i=1}^{d-6}{\cal{O}}_{\planp}(-2))
    \osum2{\cal{O}}_{\planp}(-3)\hfl{\phi}{}
(\osum_{i=1}^{d-6}{\cal{O}}_{\planp}(-1))
    \osum2{\cal{O}}_{\planp}\ra{\cal{F}}\ra 0
$$
soit de type $(II)$
$$
0\ra{\cal{O}}(-1)\osum(\osum_{i=1}^{d-5}{\cal{O}}(-2))
    \osum2{\cal{O}}_{}(-3)\hfl{\phi}{}
{\cal{O}}(-2)\osum(\osum_{i=1}^{d-5}{\cal{O}}(-1))
    \osum2{\cal{O}}\ra{\cal{F}}\ra 0.
$$
\end{prop}
La premi\`ere r\'esolution impose $d\ge 6$. Donc si $d=5$, la r\'esolution est
n\'ecessairement du deuxi\`eme type.

\begin{prop} Si $({\cal{F}},\sigma)$ est une th\^eta-caract\'eristique de
$\Theta^{1}(d)\backslash\Theta^{2}(d)$ localement libre sur son support
alors la r\'esolution minimale de
${\cal{F}}$ est de type $(I)$.
\end{prop}
\begin{proof} Il s'agit de montrer que $n_{1}=1$ est impossible dans ce
cas-l\`a. Le morphisme ${\cal{O}}_{\planp}(-1)\lra 2{\cal{O}}_{\planp}$ du
diagramme est donn\'e par deux sections $\alpha$ et $\beta$. Supposons
$\alpha$ et $\beta$ proportionnel. A ce moment l\`a, $K$ a ${\cal{O}}_{\ell}$
comme sous-faisceau,
$\ell$ la droite d\'efinie par $\alpha$. Or ce sous-faisceau s'injecte dans
${\cal{F}}$ ce qui donne une contradiction \`a la semi-stabilit\'e de
${\cal{F}}$. Sinon, les sections $\alpha$ et $\beta$ d\'efinissent deux droites
distinctes qui se coupent en un point $a\in\planp$. Soit $F$ le sous-fibr\'e
${\cal{O}}_{\planp}(-1)$ de $E$. Son $\varphi$-orthogonal est donn\'e par la
somme directe de ${\cal{O}}_{\planp}(-1)$ avec $n_{2}{\cal{O}}_{\planp}(-2)$
et le noyau de la fl\`eche $2{\cal{O}}_{\planp}(-3)\lra{\cal{O}}_{\planp}(-2)$
d\'efini par
$\alpha$ et
$\beta$. Ce noyau est isomorphe \`a ${\cal{O}}_{\planp}(-4)$. Le quotient
$E/{\cal{F}}^{\ort}$ est isomorphe \`a ${\cal{I}}_{a}(-2)$ o\`u
${\cal{I}}_{a}$ est le faisceau d'id\'eaux du point $a$. Consid\'erons le
diagramme suivant avec suites verticales  des $0$-suites, suites
horizontales exactes:
$$
\begin{diagram} &&&&0&&0\\ &&&&\sfl{}{}&&\sfl{}{}\\
0&\efl{}{}&F&\efl{}{}&E&\efl{}{}&E/F^{\ort}&\efl{}{}&0\\
&&\Vert&&\sfl{}{}&&\sfl{}{}\\
0&\efl{}{}&F&\efl{}{}&E^{\vee}&\efl{}{}&{\cal{F}}^{\vee}&\efl{}{}&0\\
&&&&\sfl{}{}&&\sfl{}{}\\ &&&&{\cal{F}}&\efl{}{}&\comp_{a}\\
&&&&\sfl{}{}&&\sfl{}{}\\ &&&&0&&0\\
\end{diagram}
$$ Chaque $0$-suite verticale n'a de la cohomologie uniquement en degr\'e $0$.
Soient ${\cal{H}}$, ${\cal{L}}$ et ${\cal{F}}'$ les faisceaux de cohomologie
respectivement. On a donc une suite exacte
$$0\lra{\cal{H}}\efl{}{}{\cal{L}}\efl{}{}{\cal{F}}'\efl{}{}0.$$ Soit
$K^{\cdot}$ le complexe
$0\ra F\ra E\ra E/F^{\ort}\ra 0$. Les suites spectrales
$'E_{2}^{p,q}=H^{p}(\ul{Ext}^{q}(K^{\cdot},\omega))$ et
$''E_{2}^{p,q}=\ul{Ext}^{p}(H^{q}(K^{\cdot},\omega))$ sont d'aboutissement
$\ul{{\tr{E}}xt}^{p+q}(K^{\cdot},{\cal{O}}_{\planp})$. On en d\'eduit la suite
exacte
$$0\efl{}{}{\cal{L}}\efl{}{}\ul{{\tr{E}}xt}^{0}(K^{\cdot},{\cal{O}}_{\planp})
\efl{}{}\comp_{a}\efl{}{}0$$ et l'identification
$\ul{{\tr{E}}xt}^{0}(K^{\cdot},{\cal{O}}_{\planp})={\cal{H}}^{\vee}.$ On
obtient un diagramme
$$
\begin{diagram} &&&&0&&0\\ &&&&\sfl{}{}&&\sfl{}{}\\
0&\efl{}{}&{\cal{H}}&\efl{}{}&{\cal{L}}&\efl{}{}&{\cal{F}}'&\efl{}{}&0\\
 &&\Vert&&\sfl{}{}&&\sfl{}{}&&\\
0&\efl{}{}&{\cal{H}}&\efl{}{}&{\cal{H}}^{\vee}&\efl{}{}&
{{\cal{E}}}&\efl{}{}&0\\ &&&&\sfl{}{}&&\sfl{}{}\\ &&&&\comp_{a}&=&\comp_{a}\\
 &&&&\sfl{}{}&&\sfl{}{}\\ &&&&0&&0\\
\end{diagram}
$$ Les th\^eta-caract\'eristiques ${\cal{F}}$ et ${{\cal{E}}}$ sont donc toutes
les deux extensions de $\comp_{a}$ avec ${\cal{F}}'$. La premi\`ere a deux, la
deuxi\`eme trois sections. Mais si ${\cal{F}}$ est localement libre sur en
$a$, l'espace de telles extensions est de dimension $1$. En particulier,
toutes ces extensions sont munies d'une forme quadratique. On obtient ainsi
une famille connexe de th\^eta-caract\'eristiques. Mais pour une telle famille
la dimension de l'espace des sections est invariante modulo
$2$, d'o\`u la contradiction cherch\'ee.
\end{proof}

\begin{cor}
La th\^eta-caract\'eristique g\'en\'erale de $\Theta^{1}\backslash\Theta^{2}$
admet
une r\'esolution de type $(I)$
\end{cor}
\begin{proof}
D'apr\`es (\cite{11} Prop. 8.3), le ferm\'e de $\Theta(d)$ des
th\^eta-caract\'eristiques non localement libres sur leur support
sch\'ematique,
est de codimension au moins $1$. Ce ferm\'e coupe, par irr\'eductibilit\'e de
$\Theta^{1}$, le ferm\'e
$\Theta^{1}$ suivant un ferm\'e de codimension au moins $2$, d'o\`u le
corollaire par la proposition pr\'ec\'edente.
\end{proof}

\np
\section{Groupe de Picard de la
composante paire}\label{GdPpair}

Consid\'erons la
composante $\Theta_{p}(d)$ des
th\^eta-caract\'eristiques semi-stables paires. On suppose, dans ce qui suit,
que $d\ge 3$.

\subsection{Le fibr\'e pfaffien}\label{defpfaff}

Soit $({\cal{F}},\sigma)$ une famille de th\^eta-caract\'eristiques
param\'etr\'ee
par la vari\'et\'e alg\'ebrique $S$. On suppose que $S$ est lisse et connexe,
que le sous-sch\'ema $S^{1}$ de $S$ (voir section \ref{BNlocus}) est de
codimension $1$ et que $S^{1}\backslash S^{2}$ est non vide. Ce sous-sch\'ema
est donc un diviseur de Weil et puisque
$S$ est lisse un diviseur de Cartier. D\'eterminons l'\'equation en un point
$s\in S^{1}\backslash S^{2}$: pour cela, on choisit une approximation
anti-sym\'etrique de la cohomologie de
${\cal{F}}$ d\'efinie au voisinage
$U_{s}$ de
$s$, qu'on peut supposer de la forme
$$\H^{0}({\cal{F}}_{s})\times U_{s}\efl{a}{}
\H^{0}({\cal{F}}_{s})^{*}\times U_{s},$$
avec $a:U_{s}\lra
\Lambda^{2}\H^{0}({\cal{F}}_{s})^{*}$.
Identifions via $\sigma$ et dualit\'e de Serre les espaces vectoriels
$\H^{0}({\cal{F}}_{s})^{*}$ avec $\H^{1}({\cal{F}}_{s})$ et consid\'erons le
morphisme
$$\Pf(a):U_{s}\lra\Lambda^{max}\H^{1}({\cal{F}}_{s})\simeq k$$
associant \`a $t\in U_{s}$ le pfaffien de la matrice anti-sym\'etrique
$a(t)$. Par d\'efinition, le sous-sch\'ema
$S^{1}$ est d\'efini au voisinage de $s\in S^{1}\backslash S^{2}$ par cette
\'equation.

Le fibr\'e inversible, unique \`a isomorphisme pr\`es, associ\'e \`a
$S^{1}$ sera appel\'e {\em fibr\'e pfaffien} et not\'e
${\cal{P}}_{({\cal{F}},\sigma)}$.
Si ${\cal{D}}_{{\cal{F}}}$ est le fibr\'e d\'eterminant associ\'e \`a la
famille
${\cal{F}}$, on a, sur l'ouvert $S^{1}\backslash S^{2}$ de $S^{1}$, que
$${\cal{P}}_{({\cal{F}},\sigma)}^{2}={\cal{D}}_{{\cal{F}}},$$ le
d\'eterminant
d'une matrice anti-sym\'etrique \'etant le carr\'e du pfaffien de celle-ci.
Par irr\'eductibilit\'e de $S^{1}$, le ferm\'e $S^{2}$ est de codimension $2$
dans
$S$. Par lissit\'e de $S$, le carr\'e du fibr\'e pfaffien est le fibr\'e
d\'eterminant
sur $S$ entier.

Consid\'erons maintenant le bon quotient
$T^{ss}_{p}(d)\lra\Theta_{p}(d)$ de la section \ref{thetaconstruction}. Soit
$({\cal{F}},\sigma,\alpha)$ la famille universelle sur $T^{ss}_{p}(d)$ et
$T^{ss,1}(d)$ le diviseur de Cartier comme ci-dessus. On verra dans la
section suivant que $\Theta^{1}(d)$ est int\`egre. Ainsi, $T^{os,1}_{p}(d)$
est irr\'eductible et comme le ferm\'e de
$T^{ss}_{p}(d)\backslash T^{os}_{p}(d)$
est de codimension $2$ dans $T^{ss}_{p}(d)$, l'ouvert
$T^{os,1}_{p}(d)\subset T^{ss,1}_{p}(d)$ est partout dense dans
$T^{ss,1}_{p}(d)$ ce qui montre que $T^{ss,1}_{p}(d)$ est irr\'eductible lui
aussi. Ce diviseur est
$GL(H)$-invariant et non-vide pour
$d\ge 5$. Soit de plus ${\cal{P}}_{({\cal{F}},\sigma)}$ le
fibr\'e pfaffien d\'efini comme ci-dessus,
${\cal{D}}_{{\cal{F}}}$ le fibr\'e d\'eterminant.
Ces fibr\'es sont des $GL(H)$-fibr\'es vectoriels inversibles.

\begin{lemme}\label{descentePD} Soit $d\ge 5$. Alors on a
\begin{list}{\arabic{lc})}{\usecounter{lc}}
\item Le fibr\'e pfaffien ${\cal{P}}_{({\cal{F}},\sigma)}$ descend \`a
$\Theta^{os}_{p}(d)$ pour
$d\geq 5$. De plus,
${\cal{P}}_{({\cal{F}},\sigma)}$ descend \`a
$\Theta_{p}(d)$ si et seulement si $d=5$.
\item Le fibr\'e d\'eterminant ${\cal{D}}_{{\cal{F}}}$ descend \`a
$\Theta_{p}(d)$.
\end{list}
\end{lemme}

On note par ${\cal{P}}$ le fibr\'e inversible obtenu par descente de
${\cal{P}}_{{\cal{F}}}$ et par
${\cal{D}}$ le fibr\'e inversible obtenu par descente de
${\cal{D}}_{{\cal{F}}}$.

\begin{proof} On utilisera le lemme de descente de
Kempf-Drezet-Narasimhan \cite{5} affirmant que si le groupe
r\'eductif $G$ op\`ere sur la vari\'et\'e alg\'ebrique
$Y$ et si $\pi:Y\lra Z$ est un bon quotient de $Y$ par $Z$, alors un
$G$-fibr\'e inversible
$\widetilde{{\cal{L}}}$ est de la forme $\pi^{*}({\cal{L}})$ pour
${\cal{L}}\in\Pic(Z)$ (\ie
$\widetilde{{\cal{L}}}$ {\em descend} \`a Z) si et seulement si l'action du
stabilisateur sur $\widetilde{{\cal{L}}}_{y}$ est triviale en tout point
$y\in Y$ d'orbite ferm\'ee.

Montrons la premi\`ere assertion: en un point
$[{\cal{F}}_{t},\sigma_{t},\alpha_{t}]\in T^{os}_{p}$ l'action du
stabilisateur $\{\pm 1\}$ sur
${\cal{P}}_{({\cal{F}}_{t},\sigma_{t})}$ est donn\'ee par
$$g\mapsto g^{\h^{1}({\cal{F}}_{t})}$$
L'action est donc triviale, $\h^{1}({\cal{F}}_{t})$ \'etant pair et, par le
lemme de Kempf, ${\cal{P}}_{({\cal{F}},\sigma)}$ descend \`a
$\Theta^{os}_{p}(d)$.
Maintenant, si $d=5$, toute th\^eta-caract\'eristique paire
$({\cal{F}},\sigma)$ de degr\'e $d$ ayant des sections est stable,
puisqu'il n'existe pas de th\^eta-caract\'eristique paire de degr\'e compris
entre $1$ et $4$ ayant des sections. L'action du stabilisateur est donc
triviale en tout point de $T^{ss}_{d}(d)$ et, par le lemme de Kempf,
${\cal{P}}_{({\cal{F}},\sigma)}$ descend \`a $\Theta_{p}(d)$. Si $d\geq 6$, on
consid\`ere la th\^eta-caract\'eristique somme directe orthogonale
$$({\cal{G}},\tau)={\cal{O}}_{C}\osum{\cal{O}}_{C'}\osum
(d-6){\cal{O}}_{\ell}$$ avec $\ell$ une droite, $C$ et $C'$ deux courbes
elliptiques non proportionnelles du plan projectif. Choisissons une
identification $\alpha:\H^{0}({\cal{G}}(N))\simeq k^{dN}$ et consid\'erons le
point $({\cal{G}},\tau,\alpha)$ de $T^{ss}_{p}(d)$. Le stabilisateur en
ce point s'identifie \`a $\{\pm 1\}\times\{\pm 1\}\times O(d-6)$ et son action
sur
$\Lambda^{2}(H^{1}({\cal{G}}))$ est donn\'ee par
$$(g_{1},g_{2},g_{3})\mapsto g_{1}^{\h^{1}({\cal{O}}_{C})}
g_{2}^{\h^{1}({\cal{O}}_{C'})}.$$
On voit alors que $(-1,1)$ et $(1,-1)$ n'op\`erent pas trivialement et que
${\cal{P}}_{({\cal{F}},\sigma)}$ ne peut descendre \'a
$\Theta^{os}_{p}(d)$ entier.
La deuxi\`eme assertion r\'esulte de \cite{9}
\end{proof}

\subsection{Le diviseur $\Theta^{1}(d)$}

Si $d\geq 5$, le diviseur $\Theta^{1}(d)$ est non-vide. Ce qui est
important pour notre propos est qu'il est int\`egre, comme le montre la
proposition suivante:

\begin{prop}
Le ferm\'e $\Theta^{1}(d)\subset\Theta_{p}(d)$ est
irr\'eductible et r\'eduit, lisse aux points ${\cal{O}}$-stables de
$\Theta^{1}(d)\backslash\Theta^{2}(d)$. \end{prop}
\begin{proof}
D\'emontrons d'abord l'\'enonc\'e concernant la lissit\'e. Si $d\leq4$,
l'\'enonc\'e est vide, supposons donc $d\geq5$. Soient
$({\cal{G}},\tau)\in\Theta^{1}(d)\backslash\Theta^{2}(d)$ une
th\^eta-caract\'eristique
${\cal{O}}$-stable, $S$ un voisinage \'etale de
$({\cal{G}},\tau)$ et $({\cal{F}},\sigma)$ une
famille universelle sur $S\times\Theta_{p}(d)$, \ie telle que le
morphisme modulaire d\'eduit de $({\cal{F}},\sigma)$ et $S$ soit le
rev\^etement
\'etale $S\lra \Theta_{p}(d)$ donn\'e. Il s'agit donc de montrer l'assertion
pour un point $s\in S$  dont l'image est
$({\cal{G}},\tau)$,
\ie tel que $({\cal{F}}_{s},\sigma_{s})\simeq({\cal{G}},\tau)$.
Choisissons, au voisinage de $s$, une approximation $E\efl{a}{}
E^{*}$ antisym\'etrique de la cohomologie de ${\cal{F}}$.

En g\'en\'eral, la vari\'et\'e d\'eterminantielle $S^{1}$ est lisse de
codimension
$1$ en $s$ si et
seulement si le morphisme canonique
$$\phi:T_{s}S\lra \Lambda^{2}(\Ker a_{s})^{*}$$
est surjectif. Comme dans \cite{6} on montre qu'on a un
diagramme anti-commutatif
$$\begin{diagram}
T_{s}S&\hfl{\phi}{}&\Lambda^{2}\H^{0}({\cal{F}}_{s})^{*}\\
\vfl{\psi}{}&\nefl{}{\chi}\\
\Ext^{1}_{asym}({\cal{F}}_{s},{\cal{F}}_{s})\\
\end{diagram}
$$
o\`u $\psi$ d\'esigne le morphisme de d\'eformation de Kodaira-Spencer et
$\chi$ le co-morphisme de Petri, \ie le morphisme d\'efini par l'accouplement
de Yoneda:
$$\H^{0}({\cal{F}}_{s})\times\Ext^{1}({\cal{F}}_{s},{\cal{F}}_{s})\lra
  \H^{1}({\cal{F}}_{s})$$

Comme $\psi$
est surjectif en $s$ il suffira de montrer que $\chi$ est surjectif en $s$.
Pour montrer que $\chi$ est surjectif il suffira de montrer que le morphisme
naturel
$$\zeta:\Ext^{1}({\cal{G}},{\cal{G}})\lra
L(\H^{0}({\cal{G}}),\H^{1}({\cal{G}}))$$ est surjectif. En effet, si $f\in
\Lambda^{2}\H^{0}({\cal{G}})^{*}$, il existerait
$g\in\Ext^{1}({\cal{G}},{\cal{G}})$ tel $\zeta(g)=f$. Or
$g'=1/2(g-^{\vee}g)\in
\Ext^{1}_{asym}({\cal{G}},{\cal{G}})$ et $\chi(g')=f$.
Pour montrer la surjectivit\'e de $\zeta$
on utilise la r\'esolution de la proposition \ref{Beilinson} de
${\cal{G}}$:
$$
0\lra F\hfl{\varphi}{} F^{\vee}\lra{\cal{G}}\lra 0.
$$
 avec  $F=(\osum_{i=1}^{d-6}{\cal{O}}_{\planp}(-2))
    \osum2{\cal{O}}_{\planp}(-3)$ ou
$F={\cal{O}}_{\planp}(-1)\osum(\osum_{i=1}^{d-6}{\cal{O}}_{\planp}(-2))
    \osum2{\cal{O}}_{\planp}(-3).$
Consid\'erons maintenant l'accouplement
$$\H^{0}(F^{\vee})\times\Ext^{2}(F^{\vee},F)\lra\H^{2}(F)$$
On obtient un diagramme commutatif
$$\begin{diagram}
\Ext^{1}({\cal{G}},{\cal{G}})&\hfl{f}{}&\Ext^{2}(F^{\vee},F)\\
\vfl{\varphi}{}&&\vfl{\vartheta}{}\\
L(\H^{0}({\cal{G}}),\H^{1}({\cal{G}}))&\hfl{g}{}&
L(\H^{0}(2{\cal{O}}_{\planp}),\H^{2}(2{\cal{O}}_{\planp}(-3))\\
\end{diagram}
$$
o\`u les fl\`eches verticales sont donn\'es par
les accouplements et o\`u $f$ est
donn\'ee par la suite spectrale qui calcule les
$\Ext^{i}({\cal{G}},{\cal{G}})$ \`a partir de la r\'esolution de ${\cal{G}}$.
Le morphisme $f$ est surjectif et $g$ est un isomorphisme. De plus, $\vartheta$
est un isomorphisme par dualit\'e de Serre. Par suite, $\varphi$ est surjectif
aussi et $S^{1}\backslash S^{2}$ est lisse aux points ferm\'es correspondant
\`a
des th\^eta-caract\'eristiques ${\cal{O}}$-stables. Ceci d\'emontre l'assertion
de
la proposition concernant la lissit\'e.

On sait d\'ej\`a que la sous-vari\'et\'e $\Theta^{1}(d)$ est irr\'eductible de
codimension $1$. D'apr\`es ce que pr\'ec\`ede ce ferm\'e est g\'en\'eriquement
lisse.
Par suite il est r\'eduit, ce qui termine la d\'emonstration de la
proposition.
\end{proof}

\subsection{L'ouvert $I(d)\subset \Theta_{p}(d)$.}
Consid\'erons le fibr\'e vectoriel
$$E=A\otimes_{k}{\cal{O}}_{\planp}(-2),$$ avec
$A=k^{d}$. Soit
$\widetilde{X}$ l'espace vectoriel des morphismes sym\'etriques
$q:E\lra E^{\vee}$. Si $\ol{q}$ est
injectif comme morphisme de faisceaux, son conoyau ${\cal{F}}$ est de
dimension $1$, muni d'une structure de th\^eta-caract\'eristique $\sigma$
provenant de la sym\'etrie de $q$. Le faisceau  ${\cal{F}}$ est forc\'ement
semi-stable. En effet, aucun sous-faisceau coh\'erent ${\cal{E}}$ de
${\cal{F}}$ ne peut avoir des sections non nulles, ce qui implique
$\chi({\cal{E}})\le 0$. On note
$X\subset\widetilde{X}$ l'ouvert des morphismes $q$ tels que $\ol{q}$ soit
injectif.

La groupe $G=GL(A)$ op\`ere sur $X$ en associant \`a au morphisme
sym\'etrique $q$ le morphisme sym\'etrique $g.q=^{t}g^{-1}\circ
q\circ g^{-1}$.

\begin{prop} Le morphisme modulaire $f:X\lra I(d)$ fait
de $I(d)$ un bon quotient de $X$ par $G$
\end{prop}
\begin{proof}
On aura besoin du lemme de transitivit\'e de bons quotients de
Seshadri:
\begin{lemme}
Soit $G$ un groupe alg\'ebrique, $N$ un sous-groupes distingu\'e ferm\'e
de $G$. Supposons que $G$ op\`ere sur les vari\'et\'es alg\'ebriques $X$ et
$Z$ et que la vari\'et\'e alg\'e\-brique $Y$ est un bon quotient de $X$ par
$N$. Soit de plus $g:X\lra Z$ un morphisme
$G$-\'equivariant
 qui se factorise suivant un morphisme $G/N$-\'equivariant
$h:Y\lra Z$. Alors $g$ est un bon quotient si et seulement si $h$
est un bon quotient.
\end{lemme}
Consid\'erons le triplet universel
$[{\cal{F}},\sigma,\alpha]$ sur
$T^{ss}_{ineff}(d)\times\planp$, o\`u $T^{ss}_{ineff}(d)$ est l'ouvert de
$T^{ss}(d)$ correspondant aux triplets dont la th\^eta-caract\'eristique
sous-jacente est ineffective. Consid\'erons de plus le fibr\'e de rep\`eres $R$
au dessus de $T^{ss}_{ineff}(d)$ dont la fibre au-dessus de
$[{\cal{F}}_{s},\sigma_{s},\alpha_{s}]$
param\`etre les isomorphismes $A\simeq\H^{1}({\cal{F}}_{s}(-1))$. Le groupe
$GL(H)\times G$ op\`ere sur
$R$ et $R\lra T^{ss}_{ineff}(d)$ est un bon quotient de $R$ par
$GL(H)$, puisque localement triviale dans la topologie de
Zariski de groupe structural $GL(H)$.
La r\'esolution \ref{Beilinson} fournit un morphisme de
$R$ dans $X$. De plus, ce morphisme fait de $R$ un fibr\'e localement
trivial dans la topologie de Zariski de groupe structural
$G$. Ces morphismes forment un diagramme commutatif
$$
\begin{diagram}
R&\lra&X\\
\vfl{}{}&&\vfl{}{}\\
T^{ss}_{ineff}(d)&\lra& I(d)
\end{diagram}
$$
En appliquant maintenant le lemme de Seshadri deux fois, obtient la
proposition.
\end{proof}

\begin{prop}\label{codim} Soit $E$ un fibr\'e vectoriel tel que
$S^{2}E^{*}(-3)$ soit engendr\'e par ses sections globales. Soit $\Gamma$
l'espace des sections globales de $S^{2}E^{*}(-3)$ et $\Gamma^{i}$ l'ouvert
des sections induisant un morphisme injectif en tant que morphisme de
faisceau $E\efl{}{}E^{\vee}$. Alors le compl\'ementaire de $\Gamma^{i}$ dans
$\Gamma$ est de codimension au moins deux.
\end{prop}
\begin{proof}
Soit $r=\rang(E)$, $\gamma=\dim \Gamma$. Consid\'erons le diagramme suivant:
$$
\begin{diagram}
\Gamma\times\planp&\efl{ev}{}&S^{2}E^{*}(-3)\\
\sfl{\pi}{}\\
\Gamma\\
\end{diagram}
$$
Le morphisme d'\'evaluation, not\'e $ev$, est un morphisme de vari\'et\'es
alg\'ebriques surjectif. Soit $D\subset S^{2}E^{*}(-3)$ le sous-fibr\'e dont la
fibre au dessus de $x\in \planp$ consiste en les $q_{x}$, tels que
$q_{x}:E_{x}\efl{}{}E^{*}_{x}$ soit de rang inf\'erieur ou \'egal \`a $r-1$.
Remarquons que $D$ est de rang $r(r+1)/2-1$. Consid\'erons l'image r\'eciproque
$\widetilde{D}$ de $D$ sous le morphisme d'\'evaluation. On a
$\widetilde{D}=\{(q,x)/\rang q(x)\le r-1\}$ et
la vari\'et\'e $\widetilde{D}$ est irr\'eductible de dimension $\gamma+1$.
Consid\'erons la projection
$\widetilde{\pi}:\widetilde{D}\lra \Gamma$
obtenue par restriction de $\pi$ \`a $\widetilde{D}$.
Il s'agit de
voir que le ferm\'e $\Gamma_{1}$ de $\Gamma$ des points o\`u la fibre est de
dimension $2$ est de codimension au moins $2$.
Remarquons d'abord que $\Gamma_{1}$ est un ferm\'e strict, la fibre
g\'en\'erale de $\widetilde{\pi}$ \'etant une courbe lisse. Maintenant, la
codimension de
$\Gamma_{1}$ dans $\Gamma$ ne peut \^etre $2$. Sinon,
$\widetilde{\pi}^{-1}(\Gamma_{1})$ serait un ferm\'e de $\widetilde{D}$ de
dimension $\gamma+1$ et serait donc, par irr\'eductibilit\'e de
$\widetilde{D}$,
\'egale \`a $\widetilde{D}$. Contradiction.
\end{proof}

Pour les estimations de codimension de cette section et de la suivante, on a
besoin de la proposition suivante, variante de la proposition 1.5 de
\cite{4}.

\begin{prop}\label{etudediff}
Soit $({\cal{F}},\sigma)$ une famille de th\^eta-caract\'eristiques
param\'etr\'ee par la vari\'et\'e alg\'ebrique $S$. Soit
$${\cal{H}}=\Groth_{isotr}({\cal{F}}/S,P)\lra S$$
le sch\'ema de Hilbert relatif des paires $(s,{\cal{E}})$ d'un point $s$
de $S$ et d'un quotient coh\'erent ${\cal{F}}_{s}/{\cal{E}}$ de polyn\^ome de
Hilbert $P$ tel que
${\cal{E}}\subset{\cal{F}}$ soit totalement isotrope pour $\sigma_{s}$.
Alors pour tout point $t=(s,{\cal{E}})$ de ${\cal{H}}$ on a la suite exacte
$$T_{t}{\cal{H}}\lra T_{s}S\efl{\omega_{+}}{}
\Ext^{1}_{asym}({\cal{E}},{\cal{E}}^{\vee}),$$
o\`u $\omega_{+}$ est le compos\'e du morphisme de Kodaira-Spencer
$\kappa:T_{s}S\lra\Ext^{1}_{asym}({\cal{F}}_{s},{\cal{F}}_{s})$ et du
morphisme canonique $\Ext^{1}_{asym}({\cal{F}}_{s},{\cal{F}}_{s})\lra
\Ext^{1}_{asym}({\cal{E}},{\cal{E}})$.
\end{prop}
\begin{lemme}\label{codim2part1}
Soit $X^{os}\subset X$ l'ouvert des $q$ induisant des
th\^eta-caract\'eristiques ${\cal{O}}$-stables. Alors la codimension
du ferm\'e compl\'ementaire \`a $X^{os}$ dans $X$ est au moins deux.
\end{lemme}
\begin{proof}
On montre le lemme en deux \'etapes: d'abord $(i)$ on montre qu'on a
$\codim_{_{X}}(X\backslash X^{s})\geq 2$, o\`u $X^{s}$ d\'esigne
l'ouvert de
$X$ correspondant aux th\^eta-caract\'eristiques stables, puis $(ii)$ on
montre que $\codim_{_{X^{s}}}(X^{s}\backslash X^{os})\geq 2$, o\`u
$X^{s}$ d\'esigne l'ouvert correspondant aux th\^eta-caract\'eristique
${\cal{O}}$-stables.
Pour $(i)$ on utilise la proposition \ref{etudediff},
appliqu\'ee \`a la famille des th\^eta-caract\'eristiques param\'etr\'ee par
$X$ et
d\'efinie par le conoyau du morphisme universel:
$$q:A\otimes{\cal{O}}_{X^{os}\times\planp}(-2)\lra
A^{*}\otimes{\cal{O}}_{X^{os}\times\planp}(-1).$$
Le morphisme $\omega_{+}$ est surjectif, la famille en question \'etant
compl\`ete et le morphisme canonique
$\Ext^{1}_{asym}({\cal{F}}_{s},{\cal{F}}_{s}^{\vee})\lra
\Ext^{1}_{asym}({\cal{E}},{\cal{E}}^{\vee})$ \'etant surjectif.

On est donc ramen\'e \'a prouver que si $0\not={\cal{E}}\subset{\cal{F}}$
est un sous-faisceau totalement isotrope de ${\cal{F}}$ et de caract\'eristique
d'Euler-Poincar\'e $0$, la dimension de
$\Ext^{1}_{asym}({\cal{E}},{\cal{E}}^{\vee})$ est au moins $2$.
Or, d'apr\`es la proposition 8.5 de $\cite{11}$, on a
$$\dim\Hom^{1}_{asym}({\cal{E}},{\cal{E}}^{\vee})-
\dim\Ext^{1}_{asym}({\cal{E}},{\cal{E}}^{\vee})=
-\frac{d({\cal{E}})(d({\cal{E}})+3)}{2}$$
d'o\`u $\dim\Ext^{1}_{asym}({\cal{E}},{\cal{E}}^{\vee})\geq 2$ pour
$d({\cal{E}})\geq 1$.

Pour $(ii)$,
consid\'erons deux sous-espaces vectoriels non nuls $A'$ et $A''$ de
$A$ de dimension respectivement $a'$ et $a''$ et tels que
$A=A'\osum A''$. On a un morphisme
$$\phi:S^{2}A^{'*}\otimes V\times S^{2}A^{''*}\otimes V
\lra S^{2}A^{*}\otimes V$$
en associant \`a $(q',q'')$ le morphisme sym\'etrique
$$\left(
\begin{matrix}
q'&0\\
0&q''\\
\end{matrix}
\right)
$$
L'image de $\phi$ est un ferm\'e de $X$ de codimension $3a'a''$. Soit
$Y\subset X$ le $G$-satur\'e de cet image. Comme l'action de $GL(A)$ sur
$X^{s}$ est propre, l'image r\'eciproque de l'image de $Y$ est \'egal \`a $Y$.
Maintenant, le sous-groupe de
$G$ des matrices du type
$$\left(
\begin{matrix} g'&0\\ 0&g''\\
\end{matrix}
\right)
$$
avec $g'\in GL(A')$ et $g''\in GL(A'')$ stabilise l'image
de $\phi$ la codimension de $Y$ est au moins
$3a'a''-a^{2}+a^{'2}+a^{''2}=a'a''.$
Pour $a=d\geq 3$ ce nombre est au moins $2$.

Il reste, \cf (\cite{11},
1.4), \`a consid\'erer le cas des th\^eta-caract\'eristiques hyperboliques.
A ce moment-l\`a, $a$ est pair et l'on consid\`ere une d\'ecomposition
$A=A'\osum A'$. Un argument analogue \`a celui ci-dessus montre que
la codimension est aussi au moins deux dans ce cas, d'o\`u $(ii)$, ce qui
termine la d\'emonstration du lemme.
\end{proof}

Soit
$\ol{G}=G/\{\pm\id\}$.
D'apr\`es ce que pr\'ec\`ede, les seules fonctions
inversibles sur $X$ sont les constantes non nulles  et le
groupe de Picard de $X$ est r\'eduit au fibr\'e en droites trivial. Il
en d\'ecoule que le groupe
$\Pic^{\ol{G}}(X)$ des  $\ol{G}$-fibr\'es inversibles sur $X$ est
isomorphe au groupe
$\chi(\ol{G})$ des caract\`eres de $\ol{G}$. Ce groupe est isomorphe \`a
$\reln$. Nous choisissons pour g\'en\'erateur le caract\`ere
$$
\chi:f\mapsto
\left\{
   \begin{diagram}\det^{-1}(f)& \text{si }d\text{ est pair}\hfill\null\\
                  \det^{-2}(f)& \text{si }d\text{ est impair}\hfill\null\\
   \end{diagram}
\right.
$$
Soit ${\cal{L}}_{\chi}$ le $\ol{G}$-fibr\'e inversible associ\'e \`a
$\chi$. Puisque $\ol{G}$ op\`ere librement sur l'ouvert $X^{os}$ des
points correspondant aux th\^eta-caract\'eristiques
${\cal{O}}$-stables, tout $\ol{G}$-fibr\'e inversible descend \`a $
I^{os}(d)$. On d\'eduit que
$\Pic( I^{os}(d))$ est isomorphe \`a $\reln$, nous choisissons
pour g\'en\'erateur le fibr\'e inversible ${\cal{L}}$, obtenu
par descente de ${\cal{L}}_{\chi}$.

\subsection{Fin de la d\'emonstration du th\'eor\`eme \protect\ref{PicPair}.}

Consid\'erons le morphisme d'oubli
$\gamma:I^{os}(d)\lra N^{s}(d,0).$

\begin{lemme}\label{PullbackL}
$\gamma^{*}({\cal{L}}_{N})=
\begin{cases}{\cal{L}}^{\otimes 2}&\text{si d est pair}\\
             {\cal{L}}&\text{si est impair}\\
\end{cases}
$
\end{lemme}
\begin{proof} On utilise la propri\'et\'e universelle qui d\'efinit
${\cal{L}}_{N}.$
Consid\'erons sur $X^{os}\times\planp$ le morphisme universel
$$q:A\otimes{\cal{O}}_{X^{os}\times\planp}(-2)\lra
A^{*}\otimes{\cal{O}}_{X^{os}\times\planp}(-1)$$
dont on note ${\cal{F}}$ le conoyau. Le
${\cal{O}}_{X^{os}\times\planp}$-module
${\cal{F}}$ d\'efinit une famille $X$-plate de faisceaux ${\cal{O}}$-stables de
dimension $1$ sur $\planp$. Soit $f_{{\cal{F}}}$ le morphisme modulaire
associ\'e \`a cette famille.
Par la propri\'et\'e universelle qui d\'efinit ${\cal{L}}_{N}$, on a
$$f^{*}_{{\cal{F}}}({\cal{L}}_{N})=\lambda_{{\cal{F}}}(u)\text{ avec $u$
la classe d'un point.}$$
Or la r\'esolution de ${\cal{F}}$ montre
que $\lambda_{{\cal{F}}}(u)$ s'identifie au $G$-fibr\'e inversible
$$(\det(A))^{-2}\otimes{\cal{O}}_{X^{os}}.$$ L'image r\'eciproque de
${\cal{L}}_{N}$ sous $f_{{\cal{F}}}$ s'identifie donc suivant si $d$ est pair
ou impair \`a ${\cal{L}}_{\chi}^{\otimes 2}$ ou ${\cal{L}}_{\chi}$.

Par construction on a un diagramme commutatif:
$$
\begin{diagram}
X^{os}\\
\sfl{\pi}{}&\sefl{}{f_{{\cal{F}}}}\\
 I^{os}(d)&\efl{\gamma}{}&N^{s}(d,0)\\
\end{diagram}
$$
Il en d\'ecoule que, suivant si $d$ est pair ou impair $\gamma^{*}$ envoie
${\cal{L}}_{N}$ sur ${\cal{L}}^{\otimes 2}$ ou ${\cal{L}}$
ce qui d\'emontre le lemme.
\end{proof}

Consid\'erons le morphisme d'oubli $\beta:\Theta^{os}(d)\lra
N^{s}(d,0)$.

\begin{lemme}\label{PullbackD} On a
$\beta^{*}({\cal{D}}_{N})=
\left\{
\begin{diagram}{\cal{O}}&
                   \text{si $d\le 4$}\hfill\null\\
               {\cal{P}}^{2}&
                   \text{si $d\ge 5$}\hfill\null\\
   \end{diagram}
\right.
$
\end{lemme}

\begin{proof}
C'est \'evident d'apr\`es la propri\'et\'e universelle qui d\'efinit le
fibr\'e d\'eterminant et le fait que ${\cal{P}}^{2}={\cal{D}}$.
\end{proof}

Consid\'erons maintenant la restriction
$$i^{*}:\Pic(\Theta^{os}_{p}(d))\lra \Pic(I^{os}(d))$$
Puisque $\Theta^{os}_{p}(d)$
est lisse, le noyau de $i^{*}$ est donn\'e par $\reln{\cal{P}}$.
Le groupe de Picard de
$\Theta_{p}^{os}(d)$  est donc librement engendr\'e par ${\cal{L}}$ et
${\cal{P}}$.

Consid\'erons maintenant le diagramme commutatif de morphismes naturels
$$\begin{diagram}
\Pic(\Theta_{p}(d))&\efl{\rho}{}&\Pic(\Theta_{p}^{os}(d))\\
\nfl{}{}& &\nfl{}{}\\
\Pic(N_{\planp}(d,0))&\efl{\rho_{N}}{}&\Pic(N_{\planp}^{s}(d,0))\\
\end{diagram}
$$
Le morphisme
$\rho_{N}$ est un
isomorphisme, $N_{\planp}(d,0)$ \'etant localement factorielle; le
morphisme $\rho$ est injectif, $\Theta_{p}(d)$ \'etant
normale. Il d\'ecoule du lemme \ref{PullbackL} et du diagramme que, si $d\ge 3$
est impair, ${\cal{L}}$ est dans l'image de $\rho$. Si $d\ge 4$ est pair
${\cal{L}}^{\otimes 2}$ est de m\^eme dans
l'image de $\rho$, mais ${\cal{L}}$ ne l'est pas puisque
${\cal{L}}_{\chi}$ ne descend pas \`a $I(d)$.
Il d\'ecoule du lemme \ref{PullbackD} et du diagramme que ${\cal{P}}^{2}$ est
dans l'image de
$\rho$ pour $d\geq 5$, mais ${\cal{P}}$ ne l'est que pour $d=5$ d'apr\`es le
lemme \ref{descentePD}. Les fibr\'e inversibles
${\cal{P}}^{2}$ et suivant si $d$ est impair ${\cal{L}}$ ou si $d$ est
pair ${\cal{L}}^{\otimes 2}$, s'\'etendent donc, et ceci de fa\c con
unique, \`a $\Theta_{p}(d)$ entier.

Comme l'image r\'eciproque de
${\cal{O}}_{\proj^{N}}(1)$ sous le morphisme sch\'ematique
$$\sigma_{_{N}}:N(d,0)\lra\proj^{M},$$ avec $M=d(d+3)$,
s'identifie \`a ${\cal{L}}_{N}$ on obtient aussi le lemme suivant:
\begin{lemme}
Si $\sigma_{\Theta}:\Theta_{p}(d)\lra\proj^{N}$, avec $N=d(d+3)/2$, est le
morphisme qui associe \`a une th\^eta-caract\'eristique son support
sch\'ematique,
l'image r\'eciproque de ${\cal{O}}_{\proj^{N}}(1)$ s'identifie, si $d$ est
pair, \`a l'extension unique de ${\cal{L}}^{2}$ \`a $\Theta_{p}(d)$, et, si $d$
est impair, \`a l'extension
 unique de ${\cal{L}}$ \`a $\Theta_{p}(d)$.
\end{lemme}

\np
\section{Groupe de Picard de la composante impaire}

Consid\'erons la
composante $\Theta_{i}(d)$ des
th\^eta-caract\'eristiques semi-stables impaires. On suppose, dans ce qui suit,
que $d\ge 4$.

\subsection{L'ouvert $ SI(d)\subset \Theta_{i}(d)$.}

Consid\'erons le fibr\'e vectoriel
$$E=H'\otimes_{k}{\cal{O}}_{\planp}(-2)\osum
   H''\otimes_{k}{\cal{O}}_{\planp}(-3),$$ avec
$H'=k^{d-3}$ et $H''=k$.
Soit de plus
$\widetilde{X}$ l'espace vectoriel des morphisme
sym\'etriques $q:E\lra E^{\vee}$.
Si $\ol{q}$ est g\'en\'eriquement injectif
son conoyau ${\cal{F}}$ est de dimension $1$,
muni d'une structure de th\^eta-caract\'eristique $\sigma$ provenant de la
sym\'etrie de $q$. On note $X'\subset\widetilde{X}$ l'ouvert des morphismes $q$
tels que $\ol{q}$ soit g\'en\'eriquement injectif.
D'apr\`es la proposition \ref{codim} la codimension du compl\'ementaire de $X'$
dans $\widetilde{X}$ est au moins $2$.
Soit $X\subset X'$ l'ouvert de $X'$ des $q$ d\'efinissant une
th\^eta-caract\'eristique semi-stable.
\begin{prop}
La codimension du ferm\'e compl\'ementaire de $X$ dans $X'$ est au moins deux.
\end{prop}
\begin{proof}
On utilise la proposition \ref{etudediff} appliqu\'ee \`a la famille de
th\^eta-caract\'eristiques
$X'$-plate d\'efinie par le quotient du morphisme universel sur
$X'\times\planp$. On v\'erifie que le morphisme $\omega_{+}$ est surjectif et
l'on est donc ramen\'e \`a prouver que si ${\cal{E}}\subset{\cal{F}}$ est un
sous-faisceau totalement isotrope non nul et propre de ${\cal{F}}$ tel que
$\chi({\cal{E}})=1$, alors $\Ext^{1}({\cal{E}},{\cal{E}}^{\vee})$ est de
dimension au moins $2$. Mais ceci se montre facilement en utilisant
une variante de la proposition 8.5 de \cite{11}.
\end{proof}

Le groupe alg\'ebrique $G=\Aut(E)$ op\`ere sur $X$ en associant \`a $g\in G$ et
$q\in X$ l'\'el\'ement $g.q=^{\vee}(g^{-1})qg^{-1}$.

\begin{prop} Le morphisme modulaire $X\lra SI(d)$ est un bon
quotient
\end{prop}

\begin{proof}
Identifions d'abord les \'el\'ements de $X$ aux matrices
$$q=\left(\begin{matrix}
u&v\\
^{\vee}v&w\\
\end{matrix}
\right)
$$
avec $u\in L(S^{2}H^{'*},V^{*})$,
$v\in \Hom(H''\otimes {\cal{O}}(-3),H^{'*}\otimes {\cal{O}}(-1))$
et $w\in L(S^{2}H^{''*},W^{*})$ avec
$V=H^{0}({\cal{O}}(1)$ et $W^{*}=\H^{0}({\cal{O}}(3))$.

Soit $N\subset G$ le
sous-groupe distingu\'e et ferm\'e de $G$ form\'e des \'el\'ements
$$g=\left(\begin{matrix}
1&b\\
0&1\\
\end{matrix}\right)
$$
avec $b\in\Hom(H''\otimes{\cal{O}}(-3),H^{\prime *}\otimes{\cal{O}}(-1))$.
Le quotient $G/N$ s'identifie \`a $GL(H')\times GL(H'')$.
Le groupe alg\'ebrique $N$ op\`ere (comme sous-groupe de $G$) sur $X$.
Cette action est donn\'ee en associant \`a $g\in N$ et $q\in X$
$$g.q=\left(\begin{matrix}
u&-ub+v
\\ -^{\vee}bu+^{\vee}v&^{\vee}bub-(^{\vee}bv+^{\vee}vb)+w\\
\end{matrix}\right)
$$
Le morphisme $u$ est g\'en\'eriquement injectif. En effet, sinon, soit $N$
son noyau. On obtiendrait un diagramme commutatif
$$\begin{diagram}
&&0&&0\\
&&\sfl{}{}&&\sfl{}{}\\
&&N&=&N\\
&&\sfl{(i,0)}{}&&\sfl{(0,^{\vee}vi)}{}\\
0&\ra&H^{'}\otimes{\cal{O}}(-2)\osum H^{''}\otimes{\cal{O}}(-3)
&\ra&H^{'*}\otimes{\cal{O}}(-1)\osum H^{''*}\otimes{\cal{O}} &\ra&
{\cal{F}}&\ra&0\\
&&\sfl{}{}&&\sfl{}{}&&\Vert\\
0&\ra&\Im(u)\osum H^{''}\otimes{\cal{O}}(-3)
&\ra&H^{'*}\otimes{\cal{O}}(-1)\osum
\frac{H^{''}\otimes{\cal{O}}}{N} &\ra&
{\cal{F}}&\ra&0\\
&&\sfl{}{}&&\sfl{}{}\\
&&0&&0\\
\end{diagram}
$$
Mais $\frac{H^{''}\otimes{\cal{O}}}{N}$ est de dimension $1$ et, comme
$\Im(u)$ est sans torsion, s'injecterait dans ${\cal{F}}$, ce qui est
impossible en raison de la semi-stabilit\'e de ${\cal{F}}$. Par cons\'equent,
si le morphisme $ub=0$, on a $b=0$. On en d\'eduit que l'action de $N$ sur
$X$ est libre, puisque l'application donn\'ee par l'op\'eration
$$\Hom(H''\otimes{\cal{O}}(-3),H^{\prime *}\otimes{\cal{O}}(-1))\times X\lra
X\times X$$
est donc une immersion ferm\'ee.

Il existe un quotient g\'eom\'etrique $Y$ de $X$ par l'action de $N$.

Consid\'erons maintenant le fibr\'e de rep\`eres $R$ au
dessus de $T^{ss}_{s-ineff}(d)$ dont la fibre au-dessus de
$[{\cal{F}}_{s},\sigma_{s},\alpha_{s}]$ param\`etre les isomorphismes
$$\beta:k^{d-3}\simeq
Ker(\H^{1}({\cal{F}}_{s}(-1))\efl{d_{1}^{-2,1}}{}
\H^{1}({\cal{F}}_{s})) \text{ et }\gamma:k\simeq\H^{1}({\cal{F}}_{s}).$$

On a un diagramme commutatif
$$\begin{diagram}
R&\efl{a}{}&X/N\\
\sfl{p}{}&&\sfl{f}{}\\
T^{ss}_{s-ineff}(d)&\efl{}{}&SI(d)
\end{diagram}
$$
o\`u le morphisme $a$ est donn\'e par la proposition \ref{Beilinson} et la
propri\'et\'e universelle \'evidente de $X/N$. Maintenant $p$ est un bon
quotient. Par le lemme de Seshadri, $R\lra SI(d)$ l'est aussi.
Puisque $a$ est un  bon quotient, $f$ l'est aussi, encore par le lemme de
Seshadri. Ce lemme applique \`a nouveau au morphisme modulaire
$X\lra SI(d)$ montre que ce dernier est un bon quotient, d'o\`u la
proposition.
\end{proof}

\begin{lemme} Soit $X^{os}\subset X$ l'ouvert des $q$ induisant des
th\^eta-caract\'eristiques ${\cal{O}}$-stables. Alors la codimension du ferm\'e
compl\'ementaire est au moins deux.
\end{lemme}
\begin{proof}
Analogue \`a la d\'emonstration du lemme \ref{codim2part1}.
\end{proof}
Soit $\ol{G}=\Aut(E)/\{\pm Id\}$. On d\'eduit des propositions
\ref{codim} et du lemme pr\'ec\'edent
que les seules fonctions inversibles sur $X$ sont les
constantes non nulles  et que le groupe de Picard de $X$ est r\'eduit au
fibr\'e
en droites trivial. Il en d\'ecoule que le groupe
$\Pic^{\ol{G}}(X)$ des  $\ol{G}$-fibr\'es inversibles sur $X$ est
isomorphe au groupe
$\chi(\ol{G})$ des caract\`eres de $\ol{G}$. Le groupe des
automorphismes de $E$  est form\'e de matrices $f$ de la forme
$$f=\left(\begin{matrix}
\alpha&\gamma \\
0 &\beta\\
\end{matrix}\right).$$ Le groupe des
caract\`eres de $\ol{G}$  est donn\'e par
$f\mapsto  (\det \alpha)^k(\det \beta)^{\ell}$ o\`u les entiers $k$  et
$\ell$ satisfont \`a l'\'equation $(d-3)k+\ell=0\bmod 2$.
Si $d$ est impair, ce sont les couples $(k,\ell)\in\reln\times\reln$
avec $\ell$ pair. Si $d$ est pair, ce sont les couples
$(k,\ell)\in\reln\times\reln$ avec $k+\ell$  pair. Nous choisissons
pour g\'en\'erateurs les caract\`eres $$
\chi_{_{1}}:f\mapsto
     \text{det}^{-1}(\alpha)\text{det}^{-(d-3)}(\beta)
$$
$$
\chi_{_{2}}:f\mapsto\text{det}^{-2}(\beta)
$$
Soient ${\cal{L}}_{1}$ et ${\cal{L}}_{2}$ les $\ol{G}$-fibr\'es
inversibles associ\'es respectivement aux caract\`eres $\chi_{_{1}}$ et
$\chi_{_{2}}$. Puisque $\ol{G}$ op\`ere librement sur l'ouvert
$X^{os}$ des points correspondant aux th\^eta-caract\'eristiques
${\cal{O}}$-stables, tout $\ol{G}$-fibr\'e inversible descend \`a
$SI^{os}(d)$. On d\'eduit que  $\Pic(SI^{os}(d))$ est isomorphe \`a un
groupe ab\'elien libre \`a deux g\'en\'erateurs. Pour g\'en\'erateurs nous
choisissons le fibr\'e inversible
${\cal{L}}$, obtenu par descente de ${\cal{L}}_{1}$ et le fibr\'e
inversible ${\cal{D}}$, obtenu par descente de ${\cal{L}}_{2}$.

\subsection{Fin de la d\'emonstration du th\'eor\`eme \ref{PicImpair}.}

Soit $d\not=6$. Alors la codimension du ferm\'e des th\^eta-caract\'eristiques
telles que $h^{0}({\cal{F}})\ge3$ est au moins deux dans $\Theta_{i}(d)$,
d'apr\`es les r\'esultats de la section $2$. Si $d=6$, ce ferm\'e est de
codimension $1$.
\begin{question} Si $d=6$, ce ferm\'e est irr\'eductible. Est-il r\'eduit?
\end{question}
On en d\'eduit que, si $d\not=6$, alors on a
$\Pic(SI^{os}(d))\simeq\Pic(\Theta_{i}^{os}(d))$,
si $d=6$, alors on a $\Pic(\Theta_{i}^{os}(d))
                \simeq\Pic(SI^{os}(d))\osum\reln$.

Consid\'erons le morphisme $\gamma:\Theta^{os}_{1}(d)\lra N_{\planp}(d,0)$.
\begin{lemme}
On a $\gamma^{*}({\cal{D}}_{N})={\cal{D}}$ et
$\gamma^{*}({\cal{L}}_{N})={\cal{L}}^{\otimes 2}
\otimes{\cal{D}}^{\otimes(d-4)}$
\end{lemme}
\begin{proof}
On utilise encore la propri\'et\'e universelle qui d\'efinit ${\cal{L}}_{N}$ et
${\cal{D}}_{N}$, cette fois-ci appliqu\'ee \`a la famille des faisceaux
semi-stables d\'efinie comme conoyau du morphisme universel
sur $X\times\planp$:
$$q:H'\otimes{\cal{O}}_{X\times\planp}(-2)\osum
   H^{\prime\prime}\otimes{\cal{O}}_{X\times\planp}(-3)\lra
   H^{\prime*}\otimes{\cal{O}}_{X\times\planp}(-1)\osum
   H^{\prime\prime*}\otimes{\cal{O}}_{X\times\planp}.$$
Si $u$ est la classe d'un point,
${\cal{L}}_{{\cal{F}}}(u)$ s'identifie \`a
$(\det(H'))^{-2}\otimes(\det(H^{\prime\prime}))^{-2}\otimes{\cal{O}}_{X}$,
si
$u={\cal{O}}_{\planp}$,  ${\cal{L}}_{{\cal{F}}}(u)$ s'identifie \`a
$(\det(H^{\prime\prime}))^{-2}\otimes{\cal{O}}_{X}$. De l\`a se d\'eduit  le
lemme, comme dans
le cas de l'\'enonc\'e analogue pour la composante paire.
\end{proof}
Consid\'erons le diagramme commutatif de morphismes naturels
$$\begin{diagram}
\Pic(\Theta_{i}(d))&\efl{\rho}{}&\Pic(\Theta_{i}^{os}(d))\\
\nfl{}{}& &\nfl{}{}\\
\Pic(N(d,0))&\efl{\rho_{_{N}}}{}&\Pic(N^{s}(d,0))\\
\end{diagram}
$$
Le morphisme
$\rho_{_{N}}:\Pic(N(d,0))\lra \Pic(N^{s}(d,0))$ est un
isomorphisme, $N(d,0)$ \'etant localement factorielle; le
morphisme $\rho$ est injectif, $\Theta(d)$ \'etant normale.
Il d\'ecoule du lemme et du diagramme que ${\cal{D}}$ est dans l'image de
$\rho_{\theta}$. De plus ${\cal{L}}^{\otimes 2}$ est dans l'image
de $\rho_{\theta}$, mais ${\cal{L}}$ ne l'est pas, puisque
${\cal{L}}_{\chi_{1}}$ ne descend pas \`a $SI(d)$ entier.
Ceci termine la d\'emonstration de th\'eor\`eme
\ref{PicImpair}.

\np
\section{D\'emonstration du th\'eor\`eme \protect\ref{UniFam}}

Soit $d\geq4$ un entier pair. Supposons l'existence d'une famille universelle
sur un ouvert $U\subset\Theta^{os}_{p}(d)$
On obtient, par image r\'eciproque, une famille universelle $({\cal{G}},\tau)$
param\'etr\'ee par l'ouvert $U'$, image r\'eciproque de $U$ sous le
morphisme
$T_{p}^{ss}(d,N)\lra\Theta_{p}(d)$.  Soit
$[{\cal{F}},\sigma,\alpha]$ le triplet universel sur
$T_{p}^{ss}(d,N)$ et consid\'erons le
faisceau ${\cal{L}}=\ul{\Hom}({\cal{G}},{\cal{F}})$. Ce faisceau est
un fibr\'e en droites sur $U'$, muni d'une action de $GL(H)$. Pour
$\alpha\in\{\pm id\}$ cette action est donn\'ee par
$$\alpha\lra\alpha\id_{{\cal{L}}}.$$ Par lissit\'e de $T_{p}^{ss}(d,N)$,
${\cal{L}}$ s'\'etend \`a
$T_{p}^{ss}(d,N)$ entier. Que l'action s'\'etende
est cons\'equence du lemme suivant, d\^u \`a Le Potier \cite{9}.

\begin{lemme}
Soit $X$ une vari\'et\'e affine, lisse et irr\'eductible sur
laquelle op\`ere le
groupe r\'eductif et connexe $G$, $f\in{\cal{O}}(X)$ un \'el\'ement non-nul,
invariant sous l'action de $G$. Soit $U_{f}$ l'ouvert d\'efini par $f\not=0$.
Alors si $L$ est un fibr\'e inversible sur $X$ toute action lin\'eaire de $G$
sur $L\restriction{U_{f}}$ s'\'etend \`a $X$ entier.
\end{lemme}

Consid\'erons maintenant la th\^eta-caract\'eristique suivante: Soit $\ell$ une
droite de $\planp$ et $({\cal{M}},\rho)$ d\'efinie par la somme directe
orthogonale
$$({\cal{M}},\rho)=\osum_{i=1}^{i=d}({\cal{O}}_{\ell}(-1),1)$$
Choisissons une identification $\beta$ de $H^{0}({\cal{M}}(N))$ avec $H$. Le
stabilisateur de $[{\cal{F}},\alpha,\sigma]$ sous l'action de $GL(H)$
s'identifie au groupe orthogonal $O(d)$. Ce stabilisateur op\`ere sur
${\cal{L}}_{[{\cal{M}},\beta,\sigma]}$. Une telle op\'eration est donn\'ee par
un caract\`ere de $O(d)$, donc de la forme
$$g\mapsto \det(g)^{n} \text{
avec } n=0,1.$$  Par densit\'e de $U'$ dans
$T_{p}^{ss}(d,N)$, l'op\'eration de $(-1)$
sur ${\cal{L}}_{[{\cal{M}},\rho,\beta]}$ est donn\'ee par $v\mapsto -v$
pour
$v\in{\cal{L}}_{[{\cal{M}},\rho,\beta]}$. Mais ceci n'est pas possible, vu
que $\det(-\id)=1$ pour $d$ pair, ce qui donne la contradiction cherch\'ee.


Soit $d\geq3$ un entier impair.

Consid\'erons la famille
universelle $[{\cal{F}},\sigma,\alpha]$ sur
$T^{os}(d,N)\times\planp$. Soient $p$ et $q$ les
projections canoniques dans le diagramme
$$
\begin{diagram}
T^{os}(d,N)\times\planp&\hfl{p_{2}}{}&\planp\\
\vfl{p_{1}}{}\\
T^{os}(d,N)\\
\end{diagram}
$$
Pour $u\in K(\planp)$, consid\'erons l'\'el\'ement suivant dans
$Pic^{GL(H)}(T^{os}(d,N)):$
$${\cal{L}}_{{\cal{F}}}(u)=\det(p_{1!}({\cal{F}}\otimes p_{2}^{*}(u)).$$
L'action de $\alpha=\pm 1$ est donn\'ee par
$$\alpha\mapsto\alpha^{<c,u>}$$
o\`u $c$ est la classe de $K_{top}(\planp)$ d\'efinie par $(0,d,0)$.
Choisissons $u=(0,d',0)$ avec $d'$ impair, ce qui donne $<c,u>=1\bmod 2$, et
posons
${\cal{L}}={\cal{L}}_{{\cal{F}}}(u)$.
On a une suite exacte
$$0\lra{\cal{A}}\lra{\cal{B}}\lra{\cal{F}}\lra 0$$
sur $T^{os}(d,N)\times\planp$ avec
${\cal{B}}=H\otimes_{_{k}}p^{*}({\cal{O}}_{\planp}(-N))$ et ${\cal{A}}$
localement libre. Consid\'erons sur $T^{os}(d,N)\times\planp$ l'injection
${\cal{A}}\otimes q^{*}({\cal{L}})\inject
{\cal{B}}\otimes q^{*}({\cal{L}})$.
Ces $Gl(H)$-fibr\'es descendent. Le conoyau sera une famille universelle
de th\^eta-caract\'eristiques  sur $\Theta^{os}(d)$,
si ${\cal{L}}$ satisfait \'a
${\cal{L}}\otimes{\cal{L}}={\cal{O}}$. Maintenant pour tout point
de $T^{os}(d,N)$ il existe
un voisinage de Zariski $U'$ tel que
${\cal{L}}\otimes{\cal{L}}$ soit trivial sur $U'$. On peut supposer cet ouvert
$GL(H)$-invariant (car ${\cal{L}}\otimes{\cal{L}}$ descend),
 d'o\`u l'existence d'une famille universelle
sur l'ouvert $U$, image de l'ouvert $U'$. Ainsi, il existe une famille
universelle de th\^eta-caract\'eristiques localement dans la topologie
de Zariski. Globalement sur $\Theta^{os}(d)$, il ne peut exister une telle
famille. Il suffit en effet de consid\'erer par exemple l'ouvert des
 th\^eta-caract\'eristiques ineffectives. D'apr\`es ce que pr\'ec\`ede
l'existence d'une telle famille supposerait l'existence d'un
$GL(A)$-fibr\'e inversible d'ordre $2$ dans $\Pic^{GL(A)}(X)$
(avec les notations de la section \ref{GdPpair}). Mais on  vu que ce
groupe est sans torsion.
\cqfd

\begin{question} Peut-on trouver un rev\^etement \'etale de degr\'e $2$ de
$\Theta(d)$ sur lequel on a une famille universelle?
\end{question}

\np


\begin{thebibliography}{99}
\bibitem{1}{\sc A. Beauville.}{\ Le groupe de monodromie des familles
universelles d'hyper\-surfaces et d'intersections compl\`etes,}{\em\  Lecture
Notes in Math. }{\bf\  1194 }{\  (1986) }{\ 8-18 }
\bibitem{2}{\sc F. Catanese.}{\ Babbage's Conjecture, Contact of Surfaces,
Symmetric Determinantal Varieties and Applications,}{\em\  Invent. math. }{\bf\
 63 }{\  (1981) }{\ 433-465 }
\bibitem{3}{\sc A. C. Dixon.}{\ Note on the reduction of a ternary quantic to a
symmetrical determinant,}{\em\  Proc. Camb. Phil. Soc. }{\bf\  11 }{\  (1902)
}{\ 350-351 }
\bibitem{4}{\sc J. M. Drezet et J. Le Potier.}{\ Fibr\'es stables et fibr\'es
exeptionnels sur le plan projectif,}{\em\  Ann. scient. Ec. Norm. Sup. }{\bf\
$4^{e}$ s\'erie, t.18 }{\  (1985) }{\ 193-244 }
\bibitem{5}{\sc J. M. Drezet et M. S. Narasimhan.}{\ Groupe de Picard des
vari\'et\'es de modules
de fibr\'es semi-stables sur les courbes alg\'ebriques.,}{\em\  Invent. math.
}{\bf\  97 }{\  (1989) }{\ 53-94 }
\bibitem{6}{\sc A. Hirschowitz.}{\ Rank techniques and jump
stratifications,}{\em\   Oxford University Press, Bombay }{\ (1987) }
\bibitem{7}{\sc Y. Laszlo.}{\ Th\'eor\`eme de Torelli g\'en\'erique pour les
intersections compl\`etes de trois quadriques de dimension paire,}{\em\
Invent. math. }{\bf\  98 }{\  (1989) }{\ 247-264 }
\bibitem{8}{\sc J. Le Potier.}{\ Fibr\'e d\'eterminant et courbes de saut sur
les surfaces alg\'ebriques,}{\em\  London Mathematical Society, Lecture Notes
Series }{\bf\  179 }{\  (1989) }{\ 213-240 }
\bibitem{9}{\sc J. Le Potier.}{\ Faisceaux semi-stables de dimension 1 sur le
plan projectif,}{\em\  preprint, Universit\'e Paris 7 }{\bf\  }{\  (1992) }
\bibitem{10}{\sc C. T. Simpson.}{\ Moduli of Representations of the Fundamental
Group of a Smooth Variety,}{\em\  Preprint, Princeton University }
\bibitem{11}{\sc C. Sorger.}{\ Th\^eta-caract\'eristiques des courbes trac\'ees
sur une
surface lisse,}{\em\ J. f\"ur die reine u. angew. Math. }{\bf\  435 }{\  (1993)
}{\
83-118 }
\bibitem{12}{\sc C. Sorger.}
{\ La semi-caract\'eristique d'Euler-Poincar\'e des faisceaux
$\omega_{_{X}}$-quadra\-tiques sur un sch\'ema de Cohen-Macaulay,}{\em\  \`a
paraitre au Bulletin de la Soci\'et\'e Math\'ematique de France }{\bf\  }{\
(1993) }
\end{thebibliography}
\end{document}